\documentclass[a4paper,tbtags,reqno,twoside]{article}
\input{isde-.sty} 

\title{Infinite-dimensional stochastic differential equations related to random matrices} 
\author{ Hirofumi Osada \\\small{
Address: Faculty of Mathematics, Kyushu University,} \\ 
\small{ Fukuoka, 819-0395, JAPAN }\\
\small{              Tel.: +81-92-802-4489 }, \quad 
\small{              Fax: +81-92-802-4489 }\\
\small{               Email: {osada@math.kyushu-u.ac.jp}  }     }

\input{isde-a.def} 
\begin{document}  \maketitle 
\begin{abstract}
We solve infinite-dimensional stochastic differential equations (ISDEs) describing an infinite number of Brownian particles interacting via two-dimensional Coulomb potentials. The equilibrium states of the associated unlabeled stochastic dynamics are the Ginibre random point field and Dyson's measures, which appear in random matrix theory. To solve the ISDEs we establish an integration by parts formula for these measures. 
Because the long-range effect of two-dimensional Coulomb potentials is 
quite strong, the properties of Brownian particles interacting 
with two-dimensional Coulomb potentials are remarkably different from those of Brownian particles interacting with Ruelle's class interaction potentials. As an example, we prove that the interacting Brownian particles associated with the Ginibre random point field satisfy plural ISDEs. 
\footnote{AMS subjectclass 2000{MSC 82C22 \and MSC 15A52 \and MSC 60J60 \and MSC 60K35 \and MSC 82B21}}
\end{abstract}

\section{ Introduction }\label{s:1}
Consider infinitely many Brownian particles $ \mathbf{X} = (X^i)_{i\in\N }$ moving in $ \Rd $ interacting via the two-dimensional (2D) Coulomb potentials 
$ \Psi _{\beta }$: 
\begin{align}\label{:11}&
\Psi _{\beta } (x) = -\beta \log |x| \quad \quad (\beta > 0)
.\end{align}
Then the stochastic dynamics 
$ \mathbf{X} = (X^i)_{i\in\N }$ is described by 
the following infinite-dimensional stochastic differential equation (ISDE):
\begin{align}\label{:12}&
dX^i_t = dB^i_t + \frac{\beta}{2} 
\limi{r} \sum_{|X_t^i-X_t^j|<r ,\, j\not=i} 
\frac{X^i_t-X^j_t}{|X^i_t-X^j_t|^{2}}dt \quad \quad (i\in\N )
.\end{align}
Here $ \{ B^i \}_{i\in\N } $ is a sequence of independent copies of $ d$-dimensional Brownian motions and $ \mathbf{X} = (X^i)_{i\in\N }$ is 
a continuous $ (\Rd )^{\N }$-valued process.

Physically this dynamics describes the motion of an infinite system of a one-component plasma in $ \Rd $. If $ d=2 $, so that the particles can be thought of as infinitely long parallel charged lines perpendicular to the confining plane \cite{forrester}. Because the Coulomb interactions $ \Psi _{\beta}$ are two-dimensional, the ISDE \eqref{:12} is meaningful only for $ d =1,2$.

The purpose of this paper is to solve the ISDE \eqref{:12} 
by relating the system to random matrix theory. 
Namely, we consider the cases $ d= 2 $, $ \beta =2$ and 
$ d=1$, $ \beta =1,2,4$. 
These are related to Ginibre ensembles ($ d=2$, $ \beta =2$) and 
Gaussian random matrices called GOE, GUE, and GSE ($ d=1$, $ \beta = 1,2,4$). 
The former is the thermodynamic limit of the distributions of eigen values of non Hermitian random Gaussian matrices, and the latter are those of orthogonal, unitary and symplectic random Gaussian matrices, respectively.

For a given interaction potential $ \Phi $, 
the study of ISDEs of this type was initiated 
by Lang \cite{lang.1}, \cite{lang.2}, and 
followed by Shiga \cite{shiga}, Fritz \cite{Fr}, Tanemura \cite{T2} and others. 
In these works $ \Phi $ is assumed to be a Ruelle's class potential, that is, 
$ \Phi $ is super stable and integrable at infinity. 
In addition, $ \Phi $ is assumed to be $ C^3_0$ 
(\cite{lang.1}, \cite{lang.2}) or to decay exponentially at infinity. 
Hence, the polynomial decay potentials have been excluded even for Ruelle's category.

We develop a new approach to solve ISDEs of this kind 
for general potentials $ \Phi $. As an application we solve \eqref{:12} 
with $ (d,\beta )$ as mentioned above. Our condition is easily checked for all Ruelle's class potentials with suitable smoothness outside the origin, so we give a new result even for this class. 

All our conditions to solve ISDEs are stated in terms of 
geometric assumptions on the ISDEs. The first step is the existence of the equilibrium state of the dynamics given by the ISDE. 
In case of Ruelle's class potentials this step is trivial because the equilibrium states are Gibbs measures, whose existence is well established in \cite{ruelle2}, and the relationship between the candidate equilibrium states and the ISDE follows from the Dobrushin-Lanford-Ruelle equation (DLR equation).  

On the other hand, when $ \Phi $ is a 2D Coulomb potential, the situation is drastically changed. Because of the unboundedness  at infinity of 2D Coulomb potentials, we can no longer use the method in \cite{ruelle2} for the construction of equilibrium states, and the DLR equation becomes meaningless. %
In the 2D Coulomb case, even the construction of infinite-volume measures for general $\beta $ has not yet been established. 
Moreover, the lack of the DLR equation requires a new device for clarifying the connection between the candidates for equilibrium states and the ISDE \eqref{:12}. For the construction, we use a result from random matrix theory \cite{mehta} and determinantal random point fields \cite{soshi.drpf}, \cite{shirai-t}. 
To clarify the relation between the measures and the ISDEs, we establish the integration by parts formula for the candidates for the equilibrium states. 
Because the candidates for the equilibrium states are given by the correlation functions defined by the determinants of some kernels, such a formula is extremely non-trivial. The calculation of such an integration by parts formula for the measures appearing in  random matrix theory is the heart of the present paper. 

The ISDE \eqref{:12} with $ d = \beta = 2 $ is the primary example of the present paper. In this case we have plural ISDEs representing the same diffusion (see Theorems \ref{l:21} and \ref{l:22}). Except for the unboundedness at infinity, the 2D Coulomb potentials have rather simple structure; they yield only repulsive force. The property of the associated stochastic dynamics is however drastically changed from that of the stochastic dynamics given by Ruelle's class potentials. Indeed, we will prove in a forth coming paper that the tagged particles are sub-diffusive. This contrasts strikingly with the result of Ruelle's class potentials \cite{o.p}. We conjecture that when $ d=2$, a phase transition occurs in $ \beta $. 

The ISDE \eqref{:12} with $ d=1$ and $ \beta = 2 $ was first investigated by Spohn \cite{sp.2}, and followed by \cite{o.dfa}, Nagao-Forrester \cite{nagao-f}, and Katori-Tanemura \cite{kt.cmp}. In these works, the dynamics was constructed by Dirichlet forms or in terms of space-time correlation functions. The ISDE was only intuitively obtained by analogy with SDEs for finite particle approximations. In this sense the ISDE \eqref{:12} has not yet been solved. 
We remark that the passage of the SDE representation from the finite particle systems to the infinite one is an extremely sensitive problem because of the long range nature of the 2D Coulomb potentials. 

It is plausible that our method is applicable to other measures appearing in random matrix theory and determinantal random point fields. We do not pursue this here. 

The organization of the paper is as follows: In \sref{s:2} we set up the mathematics and state some of the main theorems. 
In \sref{s:3} we prove Theorems \ref{l:26} and \ref{l:27}. These theorems give a general theory for solving ISDEs with long range potentials. 
In \sref{s:4} we prove \tref{l:44}, which gives a general procedure for the integration by parts formula. 
In \sref{s:5} we give a sufficient condition in \eqref{:44e}, which is a key to the integration by parts formula in \sref{s:4}. 
In \sref{s:6} we establish the integration by parts formula for the Ginibre random point field, which corresponds to the case $ d=2$ and $ \beta = 2$ in \eqref{:12}. 
In \sref{s:66} we prove Theorems \ref{l:21}--\ref{l:23}. 
In \sref{s:8} we prove the integration by parts formula for Dyson's models and complete the proof of Theorems \ref{l:24} and \ref{l:25}. These theorems correspond to the cases $ d=1 $ and $ \beta = 1,2,4$ 
in \eqref{:12}. 
In the Appendix we give the definition of the determinantal kernels of the case $ d=1$ and $ \beta = 1,4$. 

\section{ Set up and main results}\label{s:2}
Let $ \SS = \Rd $ and 
$ \SSS  = \{ \mathsf{s} = \sum _i \delta _{s_i}\, ;\, 
\mathsf{s} ( K ) < \infty \text{ for all compact sets } K \subset \SS \} $, 
where $ \delta _{a} $ stands for the delta measure at $ a $. 
We endow $ \SSS   $ with the vague topology, under which 
$ \SSS   $ is a Polish space. 
$ \SSS $ is called the configuration space over $ \SS $. 
We write $ \sss (x)=\sss (\{ x \} )$. Let 
\begin{align}\label{:20a}&
\SSSsi =
\{ \sss \in \SSS \, ;\, \, \sss (x)\le 1 \text{ for all }x \in \SS  ,\, \, 
\sss (\SS )= \infty \} 
.\end{align}
By definition, $ \SSSsi $ is the set of the configurations consisting of an infinite number of single point measures. 

For an infinite or finite product $ \SS ^k $ of $ \SS $ 
we define the map $ \ulab $ from $ \SS ^k $ to the set of measures on $ \SS $ by 
$ \ulab ((s_j)) = \sum_{j=1}^k \delta _{s_j}$. 
We omit $ k $ from the notation. 
We consider the restriction of $ \ulab $ 
on $ \ulab ^{-1} (\SSSsi )$. 
Let $ \kpath $ be the map from 
$ C([0,\infty);\SS ^k \cap \ulab ^{-1} (\SSSsi ))$ to 
$ C([0,\infty);\SSSsi )$ defined by 
\begin{align}\label{:20d}&
\kpath (\mathbf{X}) = \{ \sum_{j=1}^k \delta _{X_t^j}\}_{0\le t < \infty } 
,\end{align}
where $ \mathbf{X} = \{(X^j_t)_j\} $. We set 
$ \mathsf{X} = \kpath (\mathbf{X})$. 

A symmetric locally integrable function 
$ \map{\rho ^n }{\SS ^n}{[0,\infty ) } $ is called 
the $ n $-point correlation function of a probability measure $ \mu $ 
on $ \SSS $ w.r.t.\ the Lebesgue measure if $ \rho ^n $ satisfies 
\begin{align}\label{:28}&
\int_{A_1^{k_1}\ts \cdots \ts A_m^{k_m}} 
\rho ^n (x_1,\ldots,x_n) dx_1\cdots dx_n 
 = \int _{\SSS } \prod _{i = 1}^{m} 
\frac{\mathsf{s} (A_i) ! }
{(\mathsf{s} (A_i) - k_i )!} d\mu
 \end{align}
for any sequence of disjoint bounded measurable subsets 
$ A_1,\ldots,A_m \subset \SS $ and a sequence of natural numbers 
$ k_1,\ldots,k_m $ satisfying 
$ k_1+\cdots + k_m = n $. 
It is known 
that under a mild condition 
$ \{ \rho ^n \}_{n \in \N }$ determine 
the measure $ \mu $ \cite{soshi.drpf}.

Let $ \mug $ be the probability measure on the configuration space over $ \SS = \mathbb{R}^{2}$ 
whose $ n $-point correlation function $ \rgn $ 
w.r.t.\  the Lebesgue measure is given by 
\begin{align}\label{:21a}&
\rgn (x_1,\ldots,x_n)= 
\det [\kg (x_i,x_j)]_{1\le i,j \le n}
,\end{align} 
where $ \map{\kg }{\R ^{2}\ts \R ^{2}}{\mathbb{C}}$ 
is the kernel defined by 
\begin{align}\label{:21b}&
\kg (x,y) = \pi ^{-1} 
e^{-\frac{|x|^{2}}{2}-\frac{|y|^{2}}{2}}\cdot 
e^{x \bar{y}}
.\end{align}
Here we identify $ \R ^{2} $ as $ \mathbb{C}$ by the obvious correspondence: 
$ \R ^{2} \ni x=(x_1,x_2)\mapsto x_1 + \iii x_2 \in \mathbb{C}$, and 
$ \bar{y}=y_1-\iii y_2 $ means 
the complex conjugate under this identification, where $ \iii = \sqrt{-1}$. 
It is known that $ \mug (\SSSsi )=1 $. 
Moreover, $ \mug $ is translation and rotation invariant. 
$ \mug $ is called the Ginibre random point field. 

\begin{thm}\label{l:21} 
There exists a set $ \SSSg $ such that 
\begin{align}\label{:21c}&
\mug (\SSSg )= 1, \quad 
\SSSg \subset \SSSsi 
,\end{align}
and that, for all $  \mathbf{s}\in \ulab ^{-1}(\SSSg )$, 
there exists an $ (\R ^{2})^{\N }$-valued continuous process 
$ \mathbf{X}=(X^i)_{i\in\N }$, and 
$ (\R ^{2})^{\N }$-valued Brownian motion 
$ \mathbf{B}=(B^i)_{i\in\N }$ satisfying 
\begin{align}\label{:21d}&
dX_t^i = dB_t^i + \lim_{r\to\infty }
\sum_{|X_t^i-X_t^j|<r ,\, j\not=i}
\frac{X_t^i-X_t^j}{|X_t^i-X_t^j|^{2}} dt 
\quad (i \in\N )
,\\\label{:21e}&
\mathbf{X}_0 = \mathbf{s} 
.\end{align}
Moreover, $ \mathbf{X}=(X^i)_{i\in\N }$ satisfies 
\begin{align}\label{:21f}&
P (\mathbf{X}_t \in \ulab ^{-1}(\SSSg ),\ 0\le \forall t < \infty ) = 1 
,\\ \label{:21g} &
 P (\sup _{0\le t \le u}|X_t^{i}| < \infty \text{ for all } u, i \in \N ) = 1 
.\end{align}
\end{thm}

One specific aspect of the ISDE \eqref{:21d} is that 
its solution satisfies the second ISDE. 
Such a phenomenon never occurs in Ruelle's class potentials. 
\begin{thm}\label{l:22} 
The solution $ (\mathbf{X},\mathbf{B})$ in \tref{l:21} satisfies 
\begin{align}\label{:22a}&
dX_t^i = dB_t^i - X_t^i dt + \lim_{r\to\infty }
\sum_{|X_t^j|<r ,\, j\not=i} \frac{X_t^i-X_t^j}{|X_t^i-X_t^j|^{2}} dt 
\quad (i \in\N )
.\end{align}
\end{thm}

To clarify the meaning of the ISDEs we define the measure $ \mugone $ 
on $ \SS \ts \SSS $ by 
\begin{align}\label{:22b}&
\mugone (A\ts B) = \int_{A} \mugx (B) \rg ^{1} (x) dx 
,\end{align}
where $ \mugx = \mug (\cdot - \delta _{x}| \sss (x)\ge 1)$ 
is the Palm measure conditioned at $ x $ and $ \rg ^{1}$ 
is the 1-point correlation function of $ \mug $. 
Let $ \map{\bbb , \tilde{\bbb }}{\SS \ts \SSS }{\R ^{2} }$ be such that 
\begin{align}\label{:22c}&
\bbb (x,\mathsf{y}) 
= \limi{r} \sum_{|x-y_i|< r } \frac{x-y_i}{|x-y_i|^{2}}
,\\\label{:22d}&
\tilde{\bbb }(x,\mathsf{y}) = \limi{r} \sum_{|y_i|<r}
\frac{x-y_i}{|x-y_i|^{2}}, \qquad \text{where }\mathsf{y}=\sum_{i}\delta_{y_i}
.\end{align}
We will see in \lref{l:72} \thetag{3} that these two series converge in 
$ \Lloctwo (\mugone )$. 
We remark that neither of the series converges absolutely and, as a result, 
$ \bbb \not=\tilde{\bbb }$. 
Let 
$\mathsf{X}^{i*}_t = \sum_{j\not= i,j\in \N } \delta _{X^j_t}$. 
Then \eqref{:21d} and \eqref{:22a} can be rewritten as follows:  
\begin{align}\label{:22e}&
dX_t^i = dB_t^i + \bbb (X_t^i,\mathsf{X}^{i*}_t) dt &&(i\in\N )
,\\ \label{:22f} &
dX_t^i = dB_t^i - X_t^i dt + \tilde{\bbb }(X_t^i,\mathsf{X}^{i*}_t) dt
&&(i\in\N )
.\end{align}

A diffusion with state space $ S_0$ is a family of continuous stochastic processes with the strong Markov property starting from each point of the state space $ S_0$. So far, the unlabeled dynamics are known to be $ \SSS $-valued diffusions. We refine this as follows:

\begin{thm}	\label{l:23}
Let $ \mathbf{P}_{\mathbf{s}}$ be the distribution of 
the fully labeled dynamics $ \mathbf{X}$ given by \tref{l:21}. Then 
$\{\mathbf{P}_{\mathbf{s}}\}_{\mathbf{s}\in\ulab ^{-1}(\SSSg ) }$ 
is a diffusion with state space $ \ulab ^{-1}(\SSSg ) $.
\end{thm}

The second example is Dyson's model. 
Let $ S = \R $ and let $ \mub $ $ (\beta = 1,2,4 )$ be 
the probability measure on $ \SSS $ 
whose $ n $-point correlation function $ \rb $ 
is given by 
\begin{align} \label{:24} 
& \rb (x_1,\ldots,x_n) = \det [\Ksinb (x_i-x_j)]_{1 \le i,j \le n}
.\end{align}
Here we take $ \Ksiny (x) = {\sin ( \pi x )}/{\pi x } $. 
The definition of $ \Ksinb $ for $ \beta = 1,4 $ is given in the Appendix. 
We use quaternions to denote the kernel $ \Ksinb $ for $ \beta = 1,4 $. The precise meaning of the determinant of \eqref{:24} for $ \beta = 1,4 $ is given by \eqref{:91p}. 

The kernel $ \Ksiny $ is called the sine kernel. 
We remark that $ \Ksiny (t) = \frac{1}{2\pi }
\int_{|k| \le \pi } e ^{\iii  k t} \, dk $ 
and $ 0 \le \Ksiny \le \text{Id}$ 
as an operator on $ L^{2} (\R ) $. 

\begin{thm}\label{l:24} 
Let $ \beta =1,2,4 $. Then there exists a set 
$ \SSSdys $ such that 
\begin{align}\label{:24a}&
\mub (\SSSdys )= 1, \quad \SSSdys \subset 
\SSSsi 
,\end{align}
and that, for all $  \mathbf{s}\in \ulab ^{-1}(\SSSdys )$, 
there exists an $ \R ^{\N }$-valued continuous process 
$ \mathbf{X}=(X^i)_{i\in\N }$, and $ \R ^{\N }$-valued Brownian motion 
$ \mathbf{B}=(B^i)_{i\in\N }$ satisfying 
\begin{align} \label{:24b} 
& dX_t^i = dB_t^i + \frac{\beta }{2} \limi{r} 
\sum _{|X_t^i - X_t^j|\le r , \,  j\not = i } 
\frac{1}{X_t^i - X_t^j} dt \quad (i \in \N )
,\\\label{:24c}&
\mathbf{X}_0 = \mathbf{s} 
.\end{align}
Moreover, $ \mathbf{X}$ satisfies 
\begin{align}\label{:24f}&
P (\ulab (\mathbf{X}_t) \in \SSSdys ,\ 0\le \forall t < \infty ) = 1 
.\end{align}
\end{thm}

\begin{thm}	\label{l:25}
Let $ \SSSSdys =\ulab ^{-1}(\SSSdys )$. 
Let $ \mathbf{P}_{\mathbf{s}}$ be the distribution 
of $ \mathbf{X}$ given by \tref{l:24}. Then 
$\{\mathbf{P}_{\mathbf{s}} \}_{\mathbf{s}\in\SSSSdys }$ 
is a diffusion with state space $ \SSSSdys $.
\end{thm}


To solve the infinite-dimensional SDEs above, we prepare a general theory. 
Let $\map{\sigma }{\SoneSSS }{\R ^{d^{2}}}$ and $\map{\bbb }{\SoneSSS }{\Rd }$ 
be measurable functions. Let $ \aaa = \sigma \sigma ^{t}$. 
We assume for each $ (x,\mathsf{y})\in \SoneSSS $ 
\begin{align}\label{:26a}&
0 < \sum_{m,n=1}^{d} \aaa _{mn}(x,\mathsf{y}) \xi _{m}\xi _{n} \le 
\cref{;26x}|\xi |^{2} \quad \text{ for all } \xi =(\xi _m)
\in \Rd \backslash \{ \mathbf{0} \} 
.\end{align}
Here $ \Ct{;26x}$ is a positive constant independent of $ (x,\mathsf{y})$. 
For $ (X^i)_{i\in \N }$ we set 
$\mathsf{X}^{i*}_t = \sum_{j\not= i,j\in \N } \delta _{X^j_t}$ as before. 
Then the ISDEs we study are of the form:
\begin{align}\label{:26b}&
dX^i_t = \sigma (X^i_t,\mathsf{X}^{i*}_t) dB^i_t +
\bbb (X^i_t,\mathsf{X}^{i*}_t)dt
\quad (i\in \N )
.\end{align}
Let $ \check{\sigma }\xxxx $ be the function defined on $ \SS \ts \SS ^{\N }$ being symmetric in $ (y_j)_{j\in\N }$ for each $ x $ and satisfying  $ \check{\sigma }\xxxx =  \sigma \xxxxx $. We set $ \check{\bbb }$ similarly. Then we can rewrite \eqref{:26b} as \eqref{:26c}: 
\begin{align}\label{:26c}&
dX^i_t = \check{\sigma } (X^i_t,(X^j_t)_{j\not= i})dB^i_t + 
\check{\bbb }(X^i_t,(X^j_t)_{j\not= i})dt \quad (i\in\N )
.\end{align}
Let $ \check{\aaa }= \check{\sigma }\check{\sigma }^{t}$. Write 
$ \check{\aaa } = [\check{\aaa }_{kl}]_{1\le k,l \le d}$ and 
$ \check{\bbb } = (\check{\bbb }_{k} )_{1\le k \le d}$. 
Then intuitively the generator $ \mathbf{L}$ 
of the diffusion given by \eqref{:26c} is 
\begin{align}\label{:26d}&
\mathbf{L}= 
\frac{1}{2}\sum_{i\in\N }\sum_{k,l=1}^d
\check{\aaa }_{kl} \ssss 
\frac{\partial^{2}}{ \partial s_{ik}\partial s_{il}}
+ \sum _{i\in\N } \sum_{k=1}^d 
\check{\bbb }_k\ssss 
\frac{\partial }{ \partial s_{ik}}
.\end{align}
Here $ s_i=(s_{i1},\ldots,s_{id})\in \SS \equiv \Rd $. 

Our strategy for solving ISDE \eqref{:26b} and \eqref{:26c} 
is to use a geometric property behind the ISDE \eqref{:26b}. 
We first consider an invariant probability measure $ \mu $ 
of the unlabeled dynamics associated with \eqref{:26b}. 
Namely, we consider a probability measure $ \mu $ whose log derivative $ \dmu $ satisfies $ \bbb (x,\mathsf{y}) = \nabla _{x} \aaa (x,\mathsf{y})+
\aaa (x,\mathsf{y})\dmu (x,\mathsf{y})$. Here, to be more precise, $ \dmu $ is 
the log derivative of the measure $ \mu ^1 $ given by \eqref{:26e}, and 
the definition of $\dmu $ is given by \eqref{:26m}.  

Note that for a given pair $(\aaa , \mu )$, $ \bbb  $ is uniquely determined.  
We construct the unlabeled diffusion associated with $(\aaa , \mu )$ by using the Dirichlet space given by $(\aaa , \mu )$ and prove that the labeled process consisting of each component of the unlabeled diffusion satisfies \eqref{:26b} and \eqref{:26c}. 

If there were a Dirichlet space associated 
with the fully labeled diffusion $ \mathbf{X}=(X^i)_{i\in \N }$, 
we could use the Ito formula for each component $ X^i $ and $ X^iX^j $, 
and prove that $ \mathbf{X}$ satisfies \eqref{:26d} since 
all coordinate functions $ x^i , x^ix^j \ (i,j\in\N )$ 
would be in the domain of the Dirichlet space locally. 
We emphasize that there exist no Dirichlet spaces 
associated with the fully labeled diffusion $ \mathbf{X}$. 
Instead we introduce an infinite sequence of Dirichlet spaces 
associated with the $ k $-labeled process 
$\{((X^1_t,\ldots,X^k_t,\sum_{j>k}\delta_{X^j_t}))\}$ 
for all $ k=0,1,\ldots $. 
This sequence of $ k $-labeled processes has 
consistency and satisfies the ISDEs \eqref{:26b} and \eqref{:26c}.

Let $ \mu $ be a probability measure on 
$ (\SSS , \mathcal{B}(\SSS ))$. 
Let $ \rho ^{k}$ be the $ k $-point correlation function of 
$ \mu $ w.r.t.\ the Lebesgue measure. 
Let $ \muk $ be the measure on $ \SkS $ defined by 
\begin{align}\label{:26e}& 
\muk (A\ts B )= \int_{A}\mu _{\mathbf{x}} (B)
\rho ^{k}( \mathbf{x} ) d \mathbf{x} 
.\end{align}
Here $  \mathbf{x} =(x_1,\ldots,x_k)\in \Sk $ and 
$ d \mathbf{x} = dx_1\cdots dx_k $. 
Moreover $ \mu _{\mathbf{x}}$ is the Palm measure 
conditioned at $ \mathbf{x} =(x_1,\ldots,x_k) $ defined by 
\begin{align}\label{:26f}&
\mu _{\mathbf{x}} = \mu (\cdot - \sum_{i=1}^{k} \delta _{x_i} | \ 
\sss ( x_i )\ge 1 \text{ for }i=1,\ldots,k)
.\end{align}

We now introduce Dirichlet forms describing the $ k $-labeled dynamics. 
For a subset $ A \subset \SS $ we define the map 
$ \map{\pi _{A }}{\SSS }{\SSS } $ by 
$ \pi _{A } (\sss  ) = \sss ( A \cap \cdot )  $. 
We say a function $ \map{f}{\SSS }{\R } $ is local if $ f $ is 
$ \sigma[\pi _{ A }]$-measurable for some compact set $ A \subset \SS $. 
We say $ f $ is smooth if $ \tilde{f} $ is smooth, 
where $ \tilde{f}((s_i)) $ is the permutation invariant function in $ (s_i) $ such that 
$ f (\sss  ) = \tilde{f} ((s_i)) $ for $ \sss  = \sum _i \delta _{s_i} $. 

Let $ \di $  be the set of all local, smooth functions on $ \SSS $ with compact support. For $ f,g \in \di $ we set $ \map{\DDD [f,g]}{\SSS }{\R } $ by 
\begin{align} \label{:26g} & 
\DDD [f,g](\sss ) = 
\frac{1}{2} \sum _{ i } \sum _{m,n=1}^{d}\aaa _{mn}(s_i,\sss _i^*) 
\PD{\widetilde{f}(\mathbf{s})}{s_{im}} \PD{\widetilde{g}(\mathbf{s})}{s_{in}}
.\end{align}
Here $ \sss = \sum_{i}\delta _{s_i}$, $ \sss _i^*= \sum_{j\not=i}\delta _{s_j}$, $ s_i = (s_{i1},\ldots,s_{id}) \in \SS $, and $ \mathbf{s}=(s_i)$. 
For given $ f $ and $ g $ in $ \di $, it is easy to see that the right-hand side of \eqref{:26g} depends only on $ \sss $. So $ \DDD [f,g]$ is well defined. %
For $ f,g \in C_0^{\infty}(\Sk )\ot \di $ let $\nabla ^{\aaa ,k}[f,g]$ be the function on $ \SkS $ defined by 
\begin{align}\label{:26h}&
\nabla ^{\aaa ,k}
[f,g ] (\mathbf{x},\sss ) = 
\frac{1}{2} \sum _{j=1}^{k}
\sum_{m,n=1}^{d} \aaa _{mn} (x_j,\sum_{l\not= j}^{k}\delta_{x_l} + \sss )
\PD{f(\mathbf{x}, \sss )}{x_{jm}}\PD{g(\mathbf{x}, \sss )}{x_{jn}}
.\end{align}
where $\mathbf{x} =(x_j)\in \Sk $ and $ x_j =(x_{j1},\ldots,x_{jd})\in \SS $. 
We set $ \DDDk  $ for $ k\ge 1$ by 
\begin{align}\label{:26i}&   
\DDDk [f,g](\mathbf{x},\sss ) = 
\nabla ^{\aaa ,k}[f,g] (\mathbf{x},\sss )
+
\DDD [f(\mathbf{x},\cdot ),g(\mathbf{x},\cdot )](\sss )
.\end{align}
Let $ (\Eak ,\dnuik )$ be the bilinear form defined by  
\begin{align}\label{:26j} &
\Eak (f,g) = \int_{\SkS }  \DDDk [f,g] d\muk 
.\end{align}
When $ k=0 $, we take $ \DDDzero =\DDD $, 
$ \mu ^{0} = \mu $, and $ \Eazero = \Ea $. 
We set $ \Lm = L^2(\SSS , \mu)$ and $ \Lmuk = L^{2}(\SkS ,\muk ) $ and so on.


We assume that there exists a probability measure $ \mu $ on $ \SSS $ 
with correlation functions $ \{ \rho ^{k} \}_{k\in\N } $ satisfying 
\thetag{A.1}--\thetag{A.5}: 
\\
\thetag{A.1} $ \rho ^{k}$ is locally bounded for each $ k \in \N $. 
\\
\thetag{A.2} There exists a 
$ \dmu = (\dmu _m)_{m=1,\ldots,d}\in \{\Llocone (\mu ^{1})\}^d $ such that 
\begin{align}\label{:26m}&
\int_{\SoneSSS } \dmu f d \mu ^{1} = 
- \int_{\SoneSSS }\nabla _{x}f d \mu ^{1}
\quad \text{ for all } f \in \dione 
.\end{align}
Moreover, $ \dmu $ satisfies 
\begin{align}\label{:26n}&
\bbb  = 
\frac{1}{2}\{\nabla _{x}\aaa \} \dmu + \frac{1}{2}\aaa \dmu , 
\quad \bbb  \in \Lloctwo (\muone )
.\end{align}
Here $ \nabla _{x}f  = (\PD{f(x,\sss )}{x_m})_{m=1,\ldots,d}$ and 
$\nabla _{x}\aaa  = 
[\PD {\aaa _{mn}(x,\sss )}{ x_n}]_{m,n=1,\ldots,d }$, 
where $ x=(x_m)$. 
\\%
\thetag{A.3} 
$(\Eak ,\dnuik )$ is closable on $ \Lmuk $ for each $ k\in\{ 0 \}\cup \N $. 
\\
\thetag{A.4} 
$\mathrm{Cap}^{\mu } (\{\SSSsi \}^c) = 0 $.
\\
\thetag{A.5} 
There exists a $ T>0 $ such that for each $ R>0 $ 
\begin{align}\label{:26p}&
\liminf_{r\to \infty}\ ( \{ \int_{|x|\le {r+R}} \rho ^1 (x)dx \}
\{ \int_{\frac{r}{\sqrt{(r+R )T}}}^{\infty} e^{-u^{2}/2}du \} ) = 0
.\end{align}

\medskip

Let $ (\Eak ,\dak )$ be the closure of $ (\Eak ,\dnuik )$ on $ \Lmuk $. 
It is known \cite[Lemma 2.3]{o.tp} that $ (\Eak ,\dak )$ is quasi-regular and that the associated diffusion $ (\PPk , \mathsf{X}^k )$ exists. 
These diffusions have consistency in the sense of \eqref{:32f} and \eqref{:32g} (see \cite{o.tp}). 
We remark that $ \mathrm{Cap}^{\mu }$ in \thetag{A.4} is the capacity of the Dirichlet space $(\Eazero ,\d ^{\aaa ,0} ,\Lm )$. 
We call $ \dmu $ the log derivative of $ \mu $. 

The assumptions \thetag{A.4} and \thetag{A.5} have clear dynamical interpretations. Indeed, \thetag{A.4} means that particles never collide each other. Moreover, \thetag{A.5} means that each labeled particle never explodes \cite{o.tp}. 

\begin{thm}\label{l:26} 
Assume \thetag{A.1}--\thetag{A.5}. Then there exists an 
$ \SSS _0 $ such that 
\begin{align}\label{:26q}&
\mu (\SSS _0 )= 1, \quad \SSS _0 \subset 
\SSSsi 
,\end{align}
and that, for all $  \mathbf{s}\in \ulab ^{-1}(\SSS _0 )$, there exists 
an $\SS ^{\N }$-valued continuous process $ \mathbf{X}=(X^i)_{i\in\N }$, and 
$(\Rd )^{\N }$-valued Brownian motion $ \mathbf{B}=(B^i)_{i\in\N }$ satisfying 
\begin{align}\label{:26r}&
dX^i_t = \sigma (X^i_t,\mathsf{X}^{i*}_t)dB^i_t +
\bbb (X^i_t,\mathsf{X}^{i*}_t)dt \quad (i\in \N )
,\\ \label{:26s}&
\mathbf{X}_0 = \mathbf{s}
.\end{align}
Moreover, $ \mathbf{X}$ satisfies %
\begin{align}\label{:26t}&
P (\ulab (\mathbf{X}_t) \in \SSS_0 ,\ 0\le \forall t < \infty ) = 1 
.\end{align}
\end{thm}

\begin{rem}\label{r:26}
Let $ (\PP ^{1}, \mathsf{X}^{1})$ be the diffusion associated with 
$ (\Eaone , \daone , L ^{2}(\muone ))$. Let $ N = \{ N_{t} \} $ be the additive functional defined by $ N_t = \int_0^t \bbb  (\mathsf{X}^{1}_{u})du $. 
The assumption $ \bbb  \in \Lloctwo (\muone )$ in \eqref{:26n} is used to 
ensure that $ N = \{ N_{t} \} $ is a continuous additive functional locally of zero energy in the sense of \cite{fot}, and that $ N = \{ N_{t} \}$ possesses an increasing sequence of open sets $ \{ O_n \} $ in $ \SS \ts \SSS $ such that $ \mathrm{Cap}^{\muone }(\cup _{n}O_n^c)= 0$ and that 
\begin{align}\label{:26u}&
\limz{t} \frac{1}{t} \int_{ \SS \ts \SSS } 
\mathrm{E} ^{1}_{(x,\mathsf{y})}[N_t] \varphi (x,\mathsf{y}) d\muone 
= \int_{ \SS \ts \SSS }\bbb (x,\mathsf{y}) \varphi (x,\mathsf{y}) d\muone 
\end{align}
for any $ \varphi \in \daone $ such that $ \varphi = 0 $ on $ O_n^c $. 
Here $ \mathrm{E} ^{1}_{(x,\mathsf{y})}$ denotes the expectation w.r.t.\  the diffusion measure starting at $ (x,\mathsf{y})$. 
Indeed, the property $ \bbb  \in \Lloctwo (\muone )$ is used only for this. So we can relax the assumption that $\bbb  \in \Lloctwo (\muone )$. 
This fact will be used for Dyson's model with $ \beta = 1$ because 
$ \bbb  \in \Llocp (\muone )$ for any $ 1\le p < 2 $, but 
$ \bbb  \not\in \Lloctwo (\muone )$ in this case. 
\end{rem}

\begin{thm}	\label{l:27}
Let $ \SSSS _0 $ be the subset of $ \SS ^{\N }$ defined by 
$ \SSSS _0 =\ulab ^{-1}(\SSS _0 )$. Let $ \mathbf{P}_{\mathbf{s}}$ be the distribution of $ \mathbf{X}$ given by \tref{l:26}. Then 
$\{\mathbf{P}_{\mathbf{s}}\}_{\mathbf{s}\in\SSSS _0 }$ 
is a diffusion with state space $ \SSSS _0 $.
\end{thm}

\begin{rem}\label{r:27} \thetag{1}
There exist no Dirichlet spaces associated with the fully labeled diffusion 
$\{\mathbf{P}_{\mathbf{s}}\}_{\mathbf{s}\in\SSSS _0 }$ because the diffusion 
$\{\mathbf{P}_{\mathbf{s}}\}_{\mathbf{s}\in\SSSS _0 }$ has no invariant measures. Hence \tref{l:27} does not follow directly 
from the Dirichlet form theory. \\
\thetag{2} 
The solutions obtained in \cite{Fr}, \cite{lang.1}, \cite{lang.2}, 
\cite{shiga}, \cite{T2} for Ruelle's class interaction potentials are 
strong solutions in the sense that they are functionals of given Brownian motions. The strong Markov property of the solutions was however not proved in these works except \cite{Fr}.  It is an interesting open problem to prove that the solutions in Theorems \ref{l:21} and \ref{l:24} are strong solutions. 
\end{rem}

\begin{ex}\label{d:28} 
Let $ \Psi $ be a Ruelle's class potential, smooth outside the origin. 
Then the associated translation invariant grand canonical Gibbs measures 
constructed in \cite{ruelle2} satisfy \thetag{A.1}--\thetag{A.3} and \thetag{A.5}. \thetag{A.4} is satisfied if $ d \ge 2 $, or $ d=1$ with $ \Phi $ sufficiently repulsive at the origin \cite{inu}. More concrete examples are: \\
\thetag{1} 
Let $\Phi_{6,12}(x) = \cref{;28}\{|x|^{-12}-|x|^{-6}\}$, where $ d=3 $ 
and $\Ct{;28}>0$ is a constant. $ \Phi_{6,12}$ is called 
the Lennard-Jones 6-12 potential. The corresponding ISDE is: 
\begin{align*}&
dX^i_t = dB^i_t + \frac{\cref{;28}}{2} \sum ^{\infty}_{j=1,j\ne i} \{
\frac{12(X ^i_t-X^j_t)}{|X^i_t-X^j_t|^{14}} - 
\frac{6(X ^i_t-X^j_t)}{|X^i_t-X^j_t|^{8}}\, \} dt \quad (i\in \N )
. \end{align*}
\thetag{2} 
Let $ a > d $ and set $\Phi_a(x)=(\Ct{;29}/ a )|x|^{-a}$, where $\cref{;29}>0$. Then the corresponding ISDE is: 
\begin{align}&\label{:28b}
dX^i_t=dB^i_t + \frac{\cref{;29}}{2}
\sum ^{\infty}_{j=1,j\ne i} \frac{X ^i_t-X^j_t}{| X^i_t-X^j_t|^{a+2}}dt 
\quad (i\in \N )   
.\end{align}
At first glance the ISDE \eqref{:28b} resembles \eqref{:12} because 
 \eqref{:12} corresponds to the case $ a =0 $ in \eqref{:28b}. 
The sums in the drift terms however converge absolutely,
 unlike in \eqref{:12}. 
We emphasize that the structures of the dynamics given by the solutions of 
\eqref{:28b} and \eqref{:12} are completely different from each other. 
\end{ex}

\section{Proof of Theorems \ref{l:26} and \ref{l:27}.}\label{s:3}
In this section we prove Theorems \ref{l:26} and \ref{l:27}. 
We assume \thetag{A.1}--\thetag{A.5} throughout this section. 
Let $ (\Eak ,\dak )$ be the closure of $ (\Eak ,\dnuik )$ on $ \Lmuk $. 
We set $ \mathsf{X}^{k}= (X ^k ,\mathsf{X}) \in C([0,\infty);\Sk \ts \SSS ) $.
\begin{lem}\label{l:31}
Assume \thetag{A.1} and \thetag{A.3}. Then the following holds: \\
\thetag{1} $ (\Eak ,\dak )$ is a quasi-regular Dirichlet form on $ \Lmuk $. 
\\\thetag{2} There exists a diffusion 
$\PPk = (\{\PPk _{\xsss }\}_{\xsss \in \Sk \ts \SSS }, \mathsf{X}^{k}) $ 
associated with the Dirichlet space $\9 $. 
\end{lem}
\aaaaa 
\thetag{1} follows from Lemma 2.3 in \cite{o.tp}. 
\thetag{2} follows from \thetag{1} and Dirichlet form theory. 
\bbbbb 

Let $ \map{\lab }{\SSSsi }{\SS ^{\N }} $ be a measurable map such that 
$ \ulab \circ \lab $ is the identity map. 
We represent this map by $ \lab (\sss )= (s_1,\ldots,)$, where 
$ \sss = \sum_{i=1}^{\infty}\delta_{s_i}$. The map $ \lab $ means the label of 
the originally unlabeled particle $ \sss $ and is called a {\em label}. 
So there are infinitely many labels $ \lab $ 
satisfying the above mentioned condition. Moreover, it is easy to see that 
$ \ulab ^{-1}(\SSSsi ) = \cup _{\lab }\lab (\SSSsi )$, 
where the union is taken over all labels. 

Let $ \SSSksingle $ be the subset of $ \Sk \ts \SSS $ defined by 
$ \SSSksingle = \ulab ^{-1}(\SSSsi )$. 
For a given label $ \lab $ as above let 
$ \map{\lab _k }{\SSSsi }{\SSSksingle } $ be the map defined by 
\begin{align}\label{:32}&
\lab _k(\sum_{i=1}^{\infty}\delta_{s_i})
=(s_1,\ldots,s_k,\sum_{i=k+1}^{\infty}\delta_{s_i})
.\end{align}
Note that $ \ulab \circ \lab _k $ is the identity map. 

One can extend $ \lab $ naturally as the map from  
$ C([0,\infty);\SSSsi ) $ to $ C([0,\infty);\SS ^{\N })$. 
Indeed, for a path $\mathsf{X}=\{ \mathsf{X}_t  \}\in C([0,\infty);\SSSsi ) $, 
there exists a unique $ \{(X_t^i)\}\in C([0,\infty);\SS ^{\N })$ such that 
$ (X_0^i)=\lab (\mathsf{X}_0) $ and that 
$ \sum _i \delta _{X_t^i}=  \mathsf{X}_t $ for all $ t \in [0,\infty ) $. 
We write this map as $ \lab _{\mathrm{path}}(\mathsf{X})= \{(X_t^i)\}$. We set 
$\map{\lab _{k,\mathrm{path}}}{C([0,\infty);\SSSsi )}
{C([0,\infty);\SSSksingle )}$ similarly as 
$\lab _{\mathrm{path}}$ for $ k \ge 1$. 

We write $ \Pmt = \Pmt ^{0}  $, where $ \Pmt ^0$ is given by 
\lref{l:31}. 
\begin{lem} \label{l:32}
Assume \thetag{A.1}--\thetag{A.5}. 
Then there exists a set $ \tilde{\SSS } $ satisfying  
\begin{align}\label{:32b}&
\tilde{\SSS } \subset \SSSsi 
,\\\label{:32c}&
\mathrm{Cap}^{\mu }(\tilde{\SSS }^c)= 0 
,
\\\label{:32d}&
\Pmt (\mathsf{X}_t\in \tilde{\SSS } 
\text{ for all } t)= 1 
\quad \text{ for all } \sss \in \tilde{\SSS } 
,\\\label{:32e}&
\Pmt (\sup _{0\le t \le u}|X_t^{i}| < \infty 
\text{ for all } u, i \in \N )=1 
\quad \text{ for all } \sss \in \tilde{\SSS }
.\end{align}
Here 
$ \mathsf{X}_t= \sum_{i\in \N } \delta_{X^i_t}$. 
Moreover, for all $ k\in\N $ and any label $ \lab $ 
\begin{align}
\label{:32f} &
\PP _{\sss ^{k}}^{k}  = \PP _{\ulab (\sss ^{k})} 
\circ \lab _{k,\mathrm{path}}^{-1}
\quad \text{ for all }\sss ^{k}\in 
\lab _k (\tilde{\SSS })
,\\ \label{:32g}&
 \PP _{\sss } = \PPk _{\lab _k (\sss )}\circ \kpath ^{-1} 
\quad  \text{ for all } \sss \in \tilde{\SSS } 
.\end{align}
\end{lem}
\aaaaa 
This lemma is immediate from Theorems 2.4 and 2.5 in \cite{o.tp}. 
\bbbbb 

For $ \mathbf{s}\in \ulab ^{-1}(\tilde{\SSS } )$ 
such that $ \ulab (\mathbf{s}) = \mathsf{s}$ let 
$\mathbf{P}_{\mathbf{s}} = \mathsf{P}_{\sss }\circ \lab _{\mathrm{path}}^{-1}
$, 
where $ \lab $ is a label such that $ \lab (\mathsf{s})= \mathbf{s}$. Let 
$$ C([0,\infty);\SSSsi )_{\sss } = 
\{\mathsf{X}\in C([0,\infty);\SSSsi ); \mathsf{X}_0= \sss \}.
$$ 
We remark that $ \lab _{\mathrm{path}}|_{C([0,\infty);\SSSsi )_{\sss } } = 
\hat{\lab }_{\mathrm{path}}|_{C([0,\infty);\SSSsi )_{\sss } }$ 
for any labels $ \lab $ and $ \hat{\lab }$ satisfying 
$ \lab (\mathsf{s}) = \hat{\lab }(\mathsf{s})=\mathbf{s}$, and that 
$ \ulab ^{-1}(\tilde{\SSS }) = \cup _{\lab }\lab (\SSSsi )$. 
Hence $ \mathbf{P}_{\mathbf{s}} $ is well defined. 

\begin{lem} \label{l:33}
$ \{ \mathbf{P}_{\mathbf{s}} \}_{\mathbf{s}\in \ulab ^{-1}(\tilde{\SSS } )} $ 
is a diffusion with state space $ \ulab ^{-1}(\tilde{\SSS } )$. 
\end{lem}
\aaaaa 
We recall that $\{\PP _{\sss }\}_{\sss \in \tilde{\SSS }} $ is a diffusion 
with state space $\tilde{\SSS } $ by \lref{l:31} and \lref{l:32}. Since 
$\mathbf{P}_{\mathbf{s}}(\lab _{\mathrm{path}}(C([0,\infty);\SSSsi )_{\sss }))
= 1 $ and 
$$ \lab _{\mathrm{path}}|_{C([0,\infty);\SSSsi )_{\sss } } = 
\hat{\lab }_{\mathrm{path}}|_{C([0,\infty);\SSSsi )_{\sss } }$$
for any labels $ \lab $ and $ \hat{\lab }$ satisfying 
$ \lab (\mathsf{s}) = \hat{\lab }(\mathsf{s})=\mathbf{s}$, 
we deduce that $ \mathbf{P}_{\mathbf{s}}$ depends only on $ \PP _{\sss }$ and 
the value of the label $ \lab $ at $ \sss $. 
Hence the strong Markov property follows from that of $\{\PP _{\sss }\} $. 
The continuity of the sample paths is clear by construction.  
\bbbbb 

Let $ \aaa =[\aaa _{mn}] $ and $ \bbb $ be 
as in \eqref{:26a} and \eqref{:26b}, respectively. 
\begin{lem} \label{l:34}
Let $ M^i_t = X^i_t -X^i_0 - \int_0^t \bbb (X^i_u,\mathsf{X}^{i*}_u)du $. 
Then there exists a $ \SSS _0 \subset \tilde{\SSS }$ 
satisfying $ \mathrm{Cap}^{\mu }(\tilde{\SSS }\backslash \SSS _0 )=0 $ 
such that, for each $ \mathbf{s} \in \ulab ^{-1}(\SSS _0) $, 
the collection of the processes 
$ \{ M^i \}_{i\in\N } $ under $ \mathbf{P}_{\sss } $ 
is a sequence of $ d $-dimensional continuous local martingales such that 
\begin{align}\label{:34b}&
\langle M^i,M^j \rangle _{t} = 0 \quad (i\not=j), \quad \quad 
\langle M^i,M^i \rangle _{t} = \int_0^t 
\aaa  (X^i_u,\mathsf{X}^{i*}_u) du
.\end{align}
\end{lem}

\aaaaa 
For a diffusion process $ (P,\{ X_t \} )$ with state space $ \SS $ and 
a continuous function $ f $ on $ \SS $ we write $ A^{[f]}_t = f(X_t)-f(X_0)$. 
Then $ A^{[f]} $ becomes an additive functional (AF). An AF of this type is called a Dirichlet process. It is worthwhile to note that one can apply the Fukushima decomposition for Dirichlet processes if $ f $ is locally in the domain of the Dirichlet form associated with the diffusion. We note that $ A^{[f]}$ is defined as $ A^{[f]}_t = \tilde{f}(X_t)-\tilde{f}(X_0)$, where $ \tilde{f}$ is a quasi-continuous version of $ f $ if $ f $ is not necessary continuous but is in the domain of Dirichlet spaces. 

The process $ X^i_t-X^0_0 $ 
is an AF of the unlabeled diffusion $ (\mathsf{P}, \mathsf{X})$. 
However, $ X^i_t-X^0_0 $ is not a Dirichlet process of 
$ (\mathsf{P}, \mathsf{X})$. 
Indeed, we can not identify the position of the $ i $th particle 
without tracing all of the trajectory of the unlabeled process 
$ \mathsf{X}=\{ \mathsf{X}_u \} $ until $ u \le t $. 
On the other hand, one can regard $ X^i_t-X^0_0 $ 
as a Dirichlet process of the labeled process $ (\mathbf{P}, \mathbf{X})$ 
since the coordinate function $ x^i $ 
is a function of the state space $ \SS ^{\N }$ of $ (\mathbf{P}, \mathbf{X})$. 
However, there is no Dirichlet form associated with the labeled process 
$ (\mathbf{P}, \mathbf{X})$. 
Taking these into account, we consider the $ k $-labeled process 
$ ((X^1_t,\ldots,X^k_t, \sum_{l=k+1}^{\infty}\delta _{X^l_t}))$. 
Here $ k $ is taken such that $ i,j \le k $.
We note that the $ k $-labeled process is associated 
with the Dirichlet space $ (\Eak ,\dak , \Lmuk )$. 

Applying \cite[Theorem 5.5.1]{fot} to the function $ x^i=x^i\ot 1 \in \Rd $ 
and taking \lref{l:32} into account we deduce that there exists a set 
$ \SSS ^k_0  \subset \tilde{\SSS } $ satisfying 
$ \mathrm{Cap}^{\mu }(\tilde{\SSS }\backslash \SSS ^k_0)= 0 $ and, 
for each $ \sss \in \SSS ^k_0 $, 
the $ d $-dimensional AF $ A^{[x^i]}=\{ X^i_t-X^i_0 \} $ can be decomposed 
under $ \PPk _{\lab _k (\sss )}$ as 
\begin{align}\label{:34e}&
A^{[x^i]}= M^{[x^i]} + N^{[x^i]}
.\end{align}
Here $ M^{[x^i]}$ is a martingale AF (MAF), locally of finite energy, 
and $ N^{[x^i]}$ is a continuous AF (CAF) locally of zero energy. 
By a straightforward calculation we deduce that 
for any $ \varphi \in \dnuik $ 
\begin{align}\label{:34f}&
- \Eak (x^i,\varphi )= \int_{\SkS } 
\bbb (x^i, \sum_{l\not= i}^{k}\delta_{x^{l}} + \mathsf{y})  
\varphi (\mathbf{x} , \mathsf{y}) d \muk 
.\end{align}
Here $ \mathbf{x} = (x^{1},\ldots,x^{k}) \in \Sk $. 
By $ \bbb \in \Lloctwo (\muone )$ we see that 
$ \bbb (x^i, \sum_{l\not= i}^{k}\delta_{x^{l}} + \mathsf{y}) \in 
L^2_{\mathrm{loc}}(\muk )$. 
So, by \cite[Theorem 5.2.4]{fot} together with localization, 
we deduce that 
\begin{align}\label{:34g}&
N^{[x^i]}_t = \int_0^t \bbb (X^i_u,\sum_{l\not= i}^{k}\delta_{X^{l}_u}
 + \sum_{l=k+1}^{\infty}\delta _{X^l_u}) du
.\end{align}
Hence $ M^{[x^i]}= A^{[x^i]}-N^{[x^i]}= M^i $ 
under $ \PPk _{\lab _k (\sss )}$. This, combined with the relation 
$ (M^i,\PPk _{\lab _k (\sss )})= (M^i,\mathbf{P}_{\mathsf{s}})$ 
given by \lref{l:32}, yields that $ (M^i,\mathbf{P}_{\mathsf{s}})$ 
is a continuous local martingale. 
As for the quadratic variation of $ M^i $, we note that 
\begin{align}\label{:34h}&
\DDDk [x^i_m ,x^j_n ] (\mathbf{x},\mathsf{y}) =  
\begin{cases}0 & (i\not=j)\\
\half \aaa _{mn} (x^i,\sum_{l\not=i}^{k}\delta_{x^{l}}+ \mathsf{y}) &(i=j)
\end{cases}
.\end{align}
Here $ x^i = (x^i_m)\in \Rd $. Since 
\begin{align}\label{:34i}&
2\DDDk [x^i_m x^j_n f, x^i_m x^j_n ] - 
\DDDk [(x^i_m x^j_n)^{2} , f ] = \DDDk [x^i_m,x^j_n]\, f 
\end{align}
and $ \Eak ( f , g )=\int \DDDk [ f , g ]d\muk $, we deduce \eqref{:34b} 
from \eqref{:34h} and \cite[Theorem 5.2.3]{fot}. 

Let $ \SSS _0 = \cap _{k=1}^{\infty}\SSS ^k_0$. 
Then by \eqref{:32c} and 
$ \mathrm{Cap}^{\mu }(\tilde{\SSS }\backslash \SSS ^k_0)= 0 $ 
$ (\forall k)$ we deduce that 
$ \mathrm{Cap}^{\mu }(\tilde{\SSS }\backslash \SSS _0 )=0 $. 
Hence $\SSS _0 $ satisfies the requirement of \lref{l:34}. 
\bbbbb 

\noindent {\em Proof of \tref{l:26}}. 
For $ \mathbf{s}\in \ulab ^{-1}(\SSS  _0)$ 
let $ \mathbf{P}_{\mathbf{s}}$ as in \lref{l:34}. 
Let $ \mathbf{B}=(B^i)_{i\in\N }$ be defined by  
\begin{align}\label{:35}&
B^i_t = \int_0^t \sigma ^{-1} (X^i_u,\mathsf{X}^{i*}_u) dM^i_u 
.\end{align}
Then  $ B^i $ are $ d $-dimensional continuous local martingales. 
By \eqref{:34b} and \eqref{:35} we deduce that 
$ [\langle B^i,B^j\rangle _t]_{ i,j\in\N } = t E $. 
Here $ E $ is the unit matrix on $ (\Rd )^{\N }$. 
We deduce that $ \{ B^i \}_{i\in\N } $ are independent copies 
of $ d $-dimensional Brownian motions. Hence 
$ (\mathbf{X},\mathbf{B})$ under $ \mathbf{P}_{\mathbf{s}}$ 
is a solution of \eqref{:26r} and \eqref{:26s}. 
\eqref{:21f} follows from $ \mathrm{Cap}^{\mu }(\tilde{\SSS }\backslash \SSS _0 )=0 $. The last statement follows from \lref{l:32}, \lref{l:33} and $ \mathrm{Cap}^{\mu }(\tilde{\SSS }\backslash \SSS _0 )=0$. 
\qed

\medskip

\noindent {\em Proof of \tref{l:27}}. 
By \lref{l:33} we  see that 
$ \{ \mathbf{P}_{\mathbf{s}} \}_{\mathbf{s}\in \ulab ^{-1}(\tilde{\SSS } )} $ 
is a diffusion with state space $ \ulab ^{-1}(\tilde{\SSS } )$. 
By \lref{l:34} the set $ \SSS _0$ satisfies $ \SSS _0 \subset \tilde{\SSS }$ 
and $ \mathrm{Cap}^{\mu }(\tilde{\SSS }\backslash \SSS _0 )=0 $. 
Hence we deduce that 
$ \mathbf{P}_{\mathbf{s}}(\mathbf{X}_t\in \ulab ^{-1}(\SSS _0 )
\text{ for all }t)=1$ for each $ \mathbf{s}\in \ulab ^{-1}(\SSS _0 )$. 
So we conclude that $ \{ \mathbf{P}_{\mathbf{s}} \}_{\mathbf{s}\in \ulab ^{-1}(\SSS _0 )} $ is a diffusion with state space $ \ulab ^{-1}(\SSS  _0)$. 
\qed

\section{Log derivative of random point fields. } \label{s:4}

Let $ \mu $ be a probability measure on $ \SSS $ with locally bounded 
$ n $-point correlation function $ \rho ^{n}$ for each $ n \in \N $. 
Let $ \muone $ be the measure defined by \eqref{:26e} with $ k=1 $. 
In this section we present a sufficient condition for the existence of 
the log derivative $ \dmu $ in $ \Llocp (\muone )$ with $ 1 < p $ (\tref{l:43}) and its explicit representation (\tref{l:44}). 
We shall apply these to the Ginibre random point field and 
Dyson's model in the subsequent sections. 

We set $ \Sr = \{ x \in \SS \, ;\, |s| < \rrr \}$. 
Let $ \{ \muN \} $ be a sequence of probability measures on $ \SSS $. 
We assume that their $ n $-point correlation functions $ \{\rho ^{N,n}\} $ 
satisfy for each $ r\in \N $ 
\begin{align} \label{:41a}&
\limi{N} \rho ^{N,n} (\mathbf{x})= 
\rho ^{n} (\mathbf{x}) \quad \text{ uniformly on $\Sr ^{n}$}  
,\\\label{:41b}&
\sup_{N\in\N } \sup_{\mathbf{x}\in\Sr ^{n}} \rho ^{N,n} (\mathbf{x}) \le 
\cref{;40b} ^{-n} n ^{\cref{;40c}n}
,\end{align}
where $ 0 < \Ct{;40b}(r) < \infty $ and 
$ 0 < \Ct{;40c}(r)< 1 $ are constants independent of $ n \in \N $. 

Let $ \sigma ^{N,n}_{r}$ be 
the $ n $-density function of $ \muN $ on $ \Sr $, 
 where $ r \in \N \cup \{ \infty \} $. Then 
\begin{align}\label{:41c}& 
\sigma ^{N,n}_{r}(\mathbf{x}) = \sum_{m=0}^{\infty} \frac{(-1)^m}{m!}
\int_{\Sr ^m } \rho ^{N,n+m} (\mathbf{x},\mathbf{y}) d\mathbf{y}
.\end{align}
Let $ \sigma ^{n}_{r}$ be the $ n $-density function of $ \mu $ on $ \Sr $. 
Then the same equality as \eqref{:41c} holds. 
By \eqref{:41a}--\eqref{:41c} we deduce for each $ r\in\N $ that 
\begin{align} 
\label{:41d}&
\limi{N}\sigma ^{N,n}_{r}(\mathbf{x}) = \sigma ^{n}_{r} (\mathbf{x})
\quad \text{ uniformly on $\Sr ^{n}$ for all } n\in\N 
.\end{align}
We remark that \eqref{:41b} and \eqref{:41d} imply $\{\muN \}_{N\in\N }$ converge weakly to $\mu $. 

Let $\muNx $ be the Palm measure conditioned at $ x $ as before. 
Let $ \rho ^{N,n}_{x}$ (resp.\ $ \sigma ^{N,n}_{x,r}$) be the $ n $-point correlation (resp.\ density) function of $ \muNx $. 
Let $ \muNone $ be the measure defined by \eqref{:26e} with $ n = 1$. 
Then we deduce that 
\begin{align}\label{:41e}& 
\sigma ^{N,m}_{x,r}(\mathbf{x}) = \sum_{n=0}^{\infty} \frac{(-1)^n}{n!}
\int_{\Sr ^n } \rho ^{N,m+n}_{x} (\mathbf{x},\mathbf{y}) d\mathbf{y}
,\\\label{:41f}&
 \int f d\muNone = \sum_{n=0}^{\infty}
 \frac{1}{n!} \int_{\SS \ts \Sr ^{n}}
 \tilde{f}_n(x,\mathbf{y})
 \rho ^{N,1}(x)\sigma ^{N,n}_{x,r} (\mathbf{y})  dx d\mathbf{y} 
 .\end{align}
 Here $ f \in \dione $ and 
 $ f (x,\cdot )$ is $ \sigma [\pi _{\Sr }] $-measurable for each $ x\in \SS $. 
 Moreover, $ \tilde{f}_n(x,\mathbf{y}) $ is the function on $ \SS \ts \Sr ^{n}$ 
 being symmetric in $ \mathbf{y} = (y_1,\ldots,y_n)$ for each $ x $ and 
 $ f(x,\mathsf{y})=\tilde{f}_n(x,\mathbf{y})$ when 
 $ \mathsf{y}(\Sr ) = n $ and $ \mathsf{y} = \sum_{i=1}^n \delta _{y_i}$. 
 We set $ d\mathbf{y} = dy_1\cdots dy_n $.  It is easy to see that 
\begin{align}\label{:41g}&
\rho ^{N,n}_{x}(\mathbf{y})= \rho ^{N,1+n}(x,\mathbf{y})/\rho ^{N,1}(x) ,\quad 
\rho ^{n}_{x}(\mathbf{y})= \rho ^{1+n}(x,\mathbf{y})/\rho ^{1}(x)
.\end{align}
Here $\rho ^{n}_{x}$ is the $ n $-point correlation function of $\mu _{x} $. 
\begin{lem} \label{l:41}
Let $ \sigma ^{n}_{x,r}$ be 
the $ n $-density function of $\mu _{x} $ on $ \Sr $. 
Then for all $ n,r,s $
\begin{align} 
\label{:41h}&
\limi{N} \rho ^{N,1}(x)\sigma ^{N,n}_{x,r+s} (\mathbf{y}) = 
\rho ^{1}(x)\sigma ^{n}_{x,r+s} (\mathbf{y}) 
\quad \text{ uniformly on }\Sr \ts \SS _{r+s} ^{n} 
,\\\label{:41j}&
\limi{N}\int f d\muNone = \int f d\muone 
\quad \text{ for any }f \in C_{0}(\SS \ts \SSS ) 
,\\\label{:41k}&
\limi{n}\limsupi{N} \muNone 
(\{(x,\mathsf{y})\in \Sr \ts \SSS \, ;\, \mathsf{y}(\SS _{r+s}) \ge n\})=0
,\\\label{:41l}&
\limi{n} \muone 
(\{(x,\mathsf{y})\in \Sr \ts \SSS \, ;\, \mathsf{y}(\SS _{r+s}) \ge n\})=0
.\end{align}
\end{lem}
\aaaaa 
Combining \eqref{:41a}, \eqref{:41b}, \eqref{:41e}, and \eqref{:41g} 
implies \eqref{:41h}. \eqref{:41j} follows from \eqref{:41f} and \eqref{:41h}. 
\eqref{:41k} and \eqref{:41l} are clear by \eqref{:41a}, \eqref{:41g}, and 
the assumption that $ \rho ^{n}$ are locally bounded. 
\bbbbb 

Let $ \mathcal{B}(\SS _{\q }) $ be the Borel $ \sigma $-field of $ \SS _{\q }$. 
We regard $ \mathcal{B}(\SS _{\q }) $ as a subset of $ \mathcal{B}(\SS ) $ in 
an obvious manner and denote it by the same symbol $\mathcal{B}(\SS _{\q })$. 
Let $ \map{\varpi _{s}}{\SS \ts \SSS }{\SS \ts \SSS }$ such that 
$ \varpi _{s}(x,\mathsf{y} )=(x,\sum_{\xyi < s}\delta_{y_i})$, 
where $ \mathsf{y} = \sum_{i}\delta _{y_i}$. Let 
$$ \mathcal{F}_{\q ,s} = 
\{\mathcal{B}(\SS _{\q })\ts \mathcal{B}(\SSS )\}\cap \sigma [\varpi _{s}].
$$
Set $ \Ct{;43}(\q , N)= \muNone (\SS _{\q }\ts \SSS )$. 
Then by \eqref{:41b} $\sup_{N}\cref{;43}(\q ,N)<\infty $ for each $\q \in \N $. Without loss of generality, we can assume that $\cref{;43} > 0 $ for all $\q , N $. 
So let $ \muNonebar $ be the probability measure defined by 
$\muNonebar (\cdot )= \muNone (\cdot \cap \SS _{\q }\ts \SSS )/\cref{;43}$. 

We assume that each $\muN $ has a log derivative $\dmuN = \dmuN (x,\mathsf{y})$ such that $\dmuN - \uN \in \Llocp (\muNone )$ for some $ 1 < p < \infty $, where $ \uN = \uN (x)$ is a distribution on $\SS $. We note that $\uN $ is supposed to be independent of $\mathsf{y}\in \SSS $. Let $\bar{\dlog }^{N}_{s}\in \Llocp (\muNone )$ be such that for all $  \q \in \N $ 
\begin{align}\label{:42a}
\1 \bar{\dlog }^{N}_{s} = & \ 
\mathrm{E}^{\muNonebar }[\dmuN -\uN |\mathcal{F}_{\q ,s}]
 &&\text{for a.s.\ }\muNonebar 
\\ \notag 
= & \ 
\mathrm{E}^{\muNonebar }[\dmuN |\mathcal{F}_{\q ,s}] - \1 \uN 
&&\text{for a.s.\ }\muNonebar 
.\end{align}
Then $\{\1 \bar{\dlog }^{N}_{s}\}_{s \in\N}$ is a 
$\{\mathcal{F}_{\q ,s}\}$-martingale w.r.t.\ $ \muNonebar $ for each $ \q $. 
We remark that the second equality in \eqref{:42a} comes from the fact that 
$ \uN $ is independent of $ \mathsf{y}$. 
\begin{lem} \label{l:42} 
Let $ 1 < p < \hat{p} < \infty $. 
Assume \eqref{:41a} and \eqref{:41b}. Assume that 
\begin{align}\label{:42b}&
\cref{;42b}:= \limsup_{N\to\infty} \int_{\Sr \ts \SSS }
|\dmuN -\uN |^{\hat{p}} d\muNone < \infty 
\quad \text{ for all }\rrr \in\N 
,\end{align}
where $ \Ct{;42b}(r)$ depends only on $ r $. 
Assume that there exists a $\map{\uu }{\SS }{\Rd }$ satisfying 
\begin{align}\label{:42c}&
\limi{N}\uN =\uu \quad \text{ in } L^{\hat{p}}_{\mathrm{loc}}(\SS ,dx)
.\end{align}
Then there exists a subsequence of 
$\{\{\bar{\dlog }^{N}_{s}\}_{s \in \N }\}_{N}$, 
denoted by the same symbol, with limit $ \{ \bar{\dlog }_{s} \} _{s \in \N }$ satisfying the following:  For all $ s \in \N $ 
and $\mathcal{F}_{\q ,s}$-measurable $\varphi \in \dione $ 
\begin{align}\label{:42d}&
\int_{\Sr \ts \SSS } \bar{\dlog }_{s} \varphi \, d\muone 
= \limi{N} \int_{\Sr \ts \SSS } \bar{\dlog }^{N}_{s} \varphi \, d\muNone 
,\\ \label{:42z} &
\int_{\Sr \ts \SSS } |\bar{\dlog }_{s} |^{p} \, d\muone \le 
\liminfi{N} \int_{\Sr \ts \SSS } |\bar{\dlog }^{N}_{s}|^{p} \, d\muNone 
\le \cref{;42b}^{p/\hat{p}}
\muone (\Sr \ts \SSS )^{(\hat{p}-p)/\hat{p}} 
.\end{align}
Moreover, $ \bar{\dlog }:=\limi{s}\bar{\dlog }_{s}$ converges 
in $ \Llocp (\muone )$ and $ \muone $-almost everywhere. 
\end{lem}
\aaaaa 
By \eqref{:42a} we see that 
$ \int_{\Sr \ts \SSS } |\bar{\dlog }^{N}_{s}|^{\hat{p}} d\muNone  \le 
  \int_{\Sr \ts \SSS } |\dmuN -\uN |^{\hat{p}} d\muNone $. 
Hence by \eqref{:42b} we deduce that 
\begin{align}& \label{:42e}
\limsupi{N}\sup_{s\in\N } 
  \int_{\Sr \ts \SSS } |\bar{\dlog }^{N}_{s}|^{\hat{p}} d\muNone 
\le \cref{;42b} 
\quad \text{ for each }\q 
.\end{align}

For $ (x,\mathsf{y})\in \Sr \ts \SSS $ we write 
$ \mathsf{y}= \sum _i\delta _{y_i}$ and $ \mathbf{y}=(y_i)$. We set 
\begin{align}\label{:42f}&
\mathsf{S}_{t}^{N,m} = \{ (x,\mathsf{y})\in \Sr \ts \SSS \, ;\, 
\rho ^{N,1}(x)\sigma ^{N,m}_{x,r+s} (\mathbf{y}) < a_t 
 , \ \mathsf{y} (\SS _{r+s}) = m \} 
\end{align}
and $ \mathsf{S}_{t}^{m} $ similarly as 
$ \mathsf{S}_{t}^{N,n}$ by replacing 
$ \rho ^{N,1}(x)\sigma ^{N,n}_{x,r+s} (\mathbf{y}) $ 
by 
$ \rho ^{1}(x)\sigma ^{m}_{x,r+s} (\mathbf{y}) $. 
Here $\{ a_t \}_{t \in\N } $ is an increasing sequence 
of positive numbers such that $ \limi{t}a_{t}=\infty $ and that 
for each $ m,r,s,t, \in \N $
\begin{align}\label{:42g}&
\muone (\{ (x,\mathsf{y})\in \Sr \ts \SSS \, ;\, 
 \rho ^{1}(x)\sigma ^{m}_{x,r+s} (\mathbf{y}) = a_t \} ) = 0 
.\end{align}
We set $ \mathsf{T}_{t}^{N,n} = \bigcup_{m=1}^{n}\mathsf{S}_{t}^{N,m}$ and 
$ \mathsf{T}_{t}^{n} = \bigcup_{m=1}^{n}\mathsf{S}_{t}^{m}$. 
By \eqref{:41h}, \eqref{:42e}, \eqref{:42f}, 
and the fact that $ 1_{\Sr }\bar{\dlog }^{N}_{s}$ are 
$ \mathcal{B}(\SS _{\q })\ts \sigma [\pi _{\SS _{r+s}}]$-measurable we see that 
\begin{align}\label{:42[}&
\limi{N}|
\int_{\mathsf{T}_{t}^{N,n} } |\bar{\dlog }^{N}_{s}|^{p} d\muNone - 
\int_{\mathsf{T}_{t}^{N,n} } |\bar{\dlog }^{N}_{s}|^{p} d\muone | 
\\ \notag \le &
\limi{N}\{ \sup_{1\le m \le n}\sup_{\mathsf{S}_{t}^{N,m}}
\frac{|\rho ^{N,1}(x)\sigma ^{N,m}_{x,r+s} (\mathbf{y})  
- \rho ^{1}(x)\sigma ^{m}_{x,r+s} (\mathbf{y})|}
{\rho ^{N,1}(x)\sigma ^{N,m}_{x,r+s} (\mathbf{y}) }  \} 
\int_{\mathsf{T}_{t}^{N,n} } |\bar{\dlog }^{N}_{s}|^{p} d\muNone 
\\ \notag \le &
\limi{N}\{ \sup_{1\le m \le n}\sup_{\mathsf{S}_{t}^{N,m}}
\frac{|\rho ^{N,1}(x)\sigma ^{N,m}_{x,r+s} (\mathbf{y})  
- \rho ^{1}(x)\sigma ^{m}_{x,r+s} (\mathbf{y})|}
{a_{t}} \} 
\cref{;42b}^{p/\hat{p}} \muNone ( \Sr \ts \SSS )
^{(\hat{p}-p)/\hat{p}}\\
= & \ \ 0 
.\end{align}
By applying the H\"{o}lder inequality to $ |\bar{\dlog }^{N}_{s}|^{p} $
and by using \eqref{:42e}  we have 
\begin{align}&\label{:42h}
\int_{ \Sr \ts \SSS \backslash \mathsf{T}_{t}^{N,n} } 
|\bar{\dlog }^{N}_{s}|^{p} d\muNone \le \cref{;42b}^{p/\hat{p}} 
\muNone ( \Sr \ts \SSS \backslash \mathsf{T}_{t}^{N,n} )
^{(\hat{p}-p)/\hat{p}}
.\end{align}
By \eqref{:41h}, \eqref{:41k}, \eqref{:42g}, and $ \limi{t}a_{t}=\infty $ 
we deduce that 
\begin{align}&\label{:42i}
\limi{n}\limi{t}\limsupi{N}
\muNone ( \Sr \ts \SSS \backslash \mathsf{T}_{t}^{N,n} )
\le 
\limi{n}\limi{t}
\muone ( \Sr \ts \SSS \backslash \mathsf{T}_{t}^{n} )
=0
.\end{align}
Combining \eqref{:42h} and \eqref{:42i} we obtain 
\begin{align}\label{:42]}&
\limi{n}\limi{t}\limsupi{N}
\int_{ \Sr \ts \SSS \backslash \mathsf{T}_{t}^{N,n} } 
|\bar{\dlog }^{N}_{s}|^{p} d\muNone = 0 
.\end{align}
By \eqref{:42e}, \eqref{:42[}, and \eqref{:42]} we obtain   
\begin{align}\label{:42<}&
\limi{n}\limi{t}\limsupi{N}\int_{\mathsf{T}_{t}^{N,n} } 
|\bar{\dlog }^{N}_{s}|^{p} d\muone \le \cref{;42b}^{p/\hat{p}}
\muone (\Sr \ts \SSS )^{(\hat{p}-p)/\hat{p}} 
.\end{align}

By \eqref{:42<} we can choose a subsequence of 
$ \{ \bar{\dlog }^{N}_{s} \} $, denoted by the same symbol, 
such that $ \{1_{\mathsf{T}_{t}^{N,n}}\bar{\dlog }^{N}_{s}\}_{N \in\N}$ 
converge weakly in $ L^p(\Sr \ts \SSS ,\muone) $ to 
$ \bar{\dlog }_{s}^{t,n} $ for each $ s,t,n \in \N $. 
We can take the subsequence in such a way that the limit points 
$ \{  \bar{\dlog }_{s}^{t,n}  \} $ satisfy 
$$ 1_{\mathsf{T}_{t}^{n}}\bar{\dlog }_{s}^{t,n}=
  1_{\mathsf{T}_{t}^{n}}\bar{\dlog }_{s}^{t',n'}
\text{ for any $ n \le n'$ and $ t \le t'$. }$$
Hence we can rewrite the limit points as 
$ \bar{\dlog }_{s}^{t,n} = 1_{\mathsf{T}_{t}^{n}}\bar{\dlog }_{s}$ for any 
$ s,t,n \in \N $. By construction 
$ \limi{n}\limi{t}\bar{\dlog }_{s}^{t,n} = \bar{\dlog }_{s}$ $ \muone $-a.s.. 
Then by Fatou's lemma and \eqref{:42<} we obtain the first inequality in \eqref{:42z}. The second one is immediate from the H\"{o}lder inequality.

By \eqref{:41h} with the cut off argument similar to 
\eqref{:42[}--\eqref{:42<}  and the fact that 
$ 1_{\Sr }\bar{\dlog }^{N}_{s}\varphi $ are 
$ \mathcal{B}(\SS _{\q })\ts \sigma [\pi _{\SS _{r+s}}]$-measurable 
we obtain \eqref{:42d}. 

Let $ \bar{\mu }^{1}_{\q }=\muone (\cdot \cap \SS _{\q }\ts \SSS )/\cref{;44}$, 
where $ \Ct{;44}=\muone (\SS _{\q }\ts \SSS )$. 
By \eqref{:42a} and \eqref{:42d} we deduce that 
$ \{ \1 \bar{\dlog }_{s} \} _{s \in \N } $ is a martingale 
w.r.t.\ $ \bar{\mu }^{1}_{\q } $ for each $ \q \in \N $. 
Then the last claim follows from the martingale convergence theorems and 
\eqref{:42z}. 
\bbbbb 

\begin{thm} \label{l:43} 
Assume the same conditions as in \lref{l:42}. 
Let $ \bar{\dlog }$ be as in \lref{l:42}. 
Then the log derivative $\dmu $ of $ \mu $ exists in $ \Llocp (\muone )$ 
and is given by $ \dmu = \uu + \bar{\dlog }$. 
\end{thm}
By taking $ \uN = \uu = 0$ in \tref{l:43} we see the following:
\begin{cor}\label{l:cor} 
Let $ 1 < p < \hat{p} < \infty $. Assume \eqref{:41a} and \eqref{:41b}. 
Suppose 
\begin{align}\label{:cor}&
\limsup_{N\to\infty} \int_{\Sr \ts \SSS }|\dmuN |^{\hat{p}} d\muNone < \infty 
\quad \text{ for all }\rrr \in\N 
.\end{align}
Then the log derivative $\dmu $ of $ \mu $ exists in $ \Llocp (\muone )$. 
\end{cor}
\aaaaa 
Let $ \varphi \in \dione $. Assume, without loss of generality, that 
$\varphi $ is $\mathcal{F}_{\q ,s}$-measurable and $\varphi (x,\mathsf{y})=0 $ 
for $ x\not\in \Sq $ for some $ \q $ and $ s \in \N $. 

By \eqref{:41j} in \lref{l:41} we see that 
$\int  \nabla \varphi d\muone = \limi{N}\int  \nabla \varphi d\muNone $. 
By definition, we have 
$ - \int  \nabla \varphi \, d\muNone = \int \dmuN \varphi \, d\muNone $. 
Hence we deduce that 
\begin{align}\notag
- \int  \nabla \varphi \, d\muone & = \limi{N} \int \dmuN \varphi \, d\muNone 
\\ \notag &
= \limi{N} \int \{\uN +\bar{\dlog }^{N}_{s}\} \varphi \, d\muNone 
&&\text{ by }\eqref{:42a} 
\\ \notag &
= \int \{\uu + \bar{\dlog }_{s}\} \varphi \, d\muone 
&&\text{ by }\eqref{:42c},\ \eqref{:42d} 
\\ \notag &
= \int \{\uu +\bar{\dlog }\} \varphi \, d\muone &&\text{ by \lref{l:42}}
,\end{align}
which completes the proof.  
\bbbbb 

Let 
$\map{\g ,\, \gN ,\, \vvv ,\, \vN  }{\SS ^{2} }{\Rd }$ 
and 
$ \map{\ww }{\SS }{\Rd }$ be measurable functions. We set 
\begin{align}
\label{:44z}&
\ggs (x,\mathsf{y})= \vs + \sum_{|x-y_i|< s }\g (x,y_i)
,\\&\notag 
\ggNs (x,\mathsf{y})= \vNs + \sum_{|x-y_i|< s }\gN (x,y_i)
,\\&\notag 
\rrNs (x,\mathsf{y})= \vNsinfty + \sum_{s\le |x-y_i|}\gN (x,y_i) 
,\end{align}
where $ \mathsf{y}=\sum_{i}\delta_{y_i}$. 
We assume that 
\begin{align}\label{:44a}&
\dmuN (x,\mathsf{y})= \uN (x) + \ggNs (x,\mathsf{y}) + \rrNs (x,\mathsf{y})
,\\\label{:44b}&
\limi{N}\ggNs = \ggs \quad \text{ in }L^{\hat{p}}_{\mathrm{loc}}(\muone ) 
\quad  \text{for all }s 
,\\\label{:44e}&
\limi{s}\limsupi{N} 
  \int_{\Sr \ts \SSS } |\rrNs (x,\mathsf{y}) - \www |^{\hat{p}}  d\muNone 
= 0 , \quad  \ww \in L^{\hat{p}}_{\mathrm{loc}} (\SS ,dx)
.\end{align}

\begin{thm} \label{l:44} 
Let $ 1 < p < \hat{p} < \infty $. 
Assume \eqref{:41a}, \eqref{:41b}, and \eqref{:42c}. Assume 
\eqref{:44a}--\eqref{:44e}. 
Then the log derivative $\dmu $ exists in $ \Llocp (\muone )$ and is given by 
\begin{align}\label{:44f}&
\dmu (x,\mathsf{y})= \uu (x) + \limi{s} \ggs (x,\mathsf{y})  + \www 
.\end{align}
The convergence $ \lim \ggs $ takes place in $ \Llocp (\muone )$. 
\end{thm}
\aaaaa 
By \eqref{:44z} and \eqref{:44a} we see that $\dmuN - \uN = \ggNs + \rrNs $. 
Then \eqref{:42b} follows from \eqref{:44b} and \eqref{:44e}. 
So all the assumptions in \lref{l:42} are satisfied. 
Hence we set $ \bar{\dlog }$, $ \bar{\dlog }^{N}_{s}$ and $ \bar{\dlog }_{s}$ 
as in \lref{l:42}. 
By \tref{l:43} the log derivative 
$ \dmu = \uu + \bar{\dlog } $ exists in $\Llocp (\muone )$. 
We will prove that $ \bar{\dlog }= \limi{s} \ggs + \ww $. 

Let $\rrbar^{N}_{s} $ and $\rrbar^{N}_{ss} $ be 
functions such that for all $\q $ 
\begin{align}\label{:44g}&
\1 \rrbar^{N}_{s} = \1 \mathrm{E}^{\muNonebar }[\rrN_{0} |\mathcal{F}_{\q ,s}], \quad 
\1 \rrbar^{N}_{ss} = \1 \mathrm{E}^{\muNonebar }[\rrN_{s} |\mathcal{F}_{\q ,s}]
.\end{align}
Then $\{\1 \rrbar^{N}_{s}\}_{s \in\N}$ is a martingale w.r.t.\ $ \muNonebar $ for all $ \q $. By the second equality in \eqref{:44g} combined with \eqref{:42a} and $\dmuN =\uN + \ggNs + \rrNs $ we obtain 
\begin{align}\label{:44h}&
 \bar{\dlog }^{N}_{s} - \ggNs - \ww = \rrbar^{N}_{ss} - \ww 
.\end{align}

By \lref{l:42} we see that $\{\bar{\dlog }^{N}_{s}\}_{N}$ converge weakly to 
$\bar{\dlog }^{N}_{s}$ in $ \Llocp (\muone )$. Hence we deduce that 
\begin{align}\label{:44i}
\int_{\Sr \ts \SSS } 
&| \bar{\dlog }_{s} - \ggs - \ww |^{p}  d\muone 
\le 
\liminf_{N\to\infty} \int_{\Sr \ts \SSS } 
|\bar{\dlog }^{N}_{s} - \ggs - \ww |^{p} d\muone 
.\end{align}
By the cut off argument similar to \eqref{:42[}--\eqref{:42<} we deduce that 
\begin{align}\label{:42k}&\quad 
\limsup_{N\to\infty} \int_{\Sr \ts \SSS } 
|\bar{\dlog }^{N}_{s} - \ggs - \ww |^{p} d\muone 
\\ \notag &
\le \limi{n}\limi{t}\limsupi{N}
\int_{\mathsf{T}_{t}^{N,n} } 
|\bar{\dlog }^{N}_{s} - \ggs - \ww |^{p} d\muNone 
\\ \notag &
\le \limsupi{N}
\int_{\Sr \ts \SSS } 
|\bar{\dlog }^{N}_{s} - \ggs - \ww |^{p} d\muNone 
&&\text{ by }
\mathsf{T}_{t}^{N,n} \subset \Sr \ts \SSS 
\\ \notag &
= 
\limsup_{N\to\infty} \int_{\Sr \ts \SSS } 
|\bar{\dlog }^{N}_{s} - \ggNs - \ww |^{p} d\muNone 
&&\text{ by }\eqref{:44b} 
\\ \notag & = 
\limsup_{N\to\infty} \int_{\Sr \ts \SSS } 
|\rrbar^{N}_{ss} - \ww |^{p} d\muNone 
&&\text{ by }\eqref{:44h} 
.\end{align}
By \eqref{:44e} and \eqref{:44g} together with the assumption 
that $ \ww = \www $ is independent of $ \mathsf{y}$ we easily see that 
\begin{align}\label{:44j}&
\limi{s}\limsup_{N\to\infty} \int_{\Sr \ts \SSS } 
|\rrbar^{N}_{ss} - \ww |^{p} d\muNone  = 0  
.\end{align}

Putting \eqref{:44i}, \eqref{:42k} and \eqref{:44j} together we obtain 
\begin{align*}&
\limi{s} \int_{\Sr \ts \SSS } 
| \bar{\dlog }_{s} - \ggs - \ww |^{p}  d\muone = 0  
.\end{align*}
This combined with 
$ \dmu = \uu + \bar{\dlog } =\uu + \lim_{s}\bar{\dlog }_{s}$ 
implies \eqref{:44f}. Because the convergence of $ \lim \bar{\dlog }_{s}$ takes place in $ \Llocp (\muone )$, so does the convergence of $ \lim \ggs $. 
\bbbbb 

\section{Sufficient conditions for \eqref{:44e}} \label{s:5}
The purpose of this section is to give sufficient conditions for 
$ \ggs \in L^2_{\mathrm{loc}}(\muone ) $ and \eqref{:44e} in terms of correlation functions. 
\begin{lem} \label{l:51} 
Let $ \g _{s} (x,y) = 1_{\Ss }(x-y) \g (x,y) $. 
Then $ \ggs \in L^2_{\mathrm{loc}}(\muone ) $ follows from 
\begin{align}\label{:51a}&
\int_{\Sr }|\vs |^2\rho ^{1}(x)dx + 
\int_{\Sr \ts \SS } |\g _{s}(x,y)|^{2} \rho ^{2}(x,y)dxdy 
\\ \notag &
+ \int_{\Sr \ts \S^{2} }\g _{s}(x,y)\cdot \g _{s}(x,z) \rho ^{3}(x,y,z)dxdydz 
<\infty 
\quad \text{ for all } r  \in \N 
.\end{align}
Here $ \cdot $ denotes the standard inner product of $ \Rd $. 
\end{lem}
\aaaaa 
By definition $ \ggs (x,\mathsf{y})= \vs + \sum_{|x-y_i|< s }\g (x,y_i)$. 
By $ d\muone = \rho ^{1}(x)\mu _{x} dx $, \eqref{:41g}, 
and a simple calculation of correlation functions we see that 
\begin{align*}&
 \int_{\Sr \ts \SSS } |\sum_{|x-y_i|< s }\g (x,y_i)|^{2} d\muone = 
\int_{\Sr }
\mathrm{E}^{\mu _{x} } [|\sum_{|x-y_i|< s }\g (x,y_i)|^{2}]
\rho ^{1}(x)dx 
\\ \notag & = \int_{\Sr }
\{ \int_{\SS ^{2}}\g _{s}(x,y)\cdot \g _{s}(x,z) \rho ^{2}_{x}(y,z) dydz + 
\int_{\SS }|\g _{s}(x,y)|^{2} \rho ^{1}_{x}(y)dy \} \rho ^{1}(x)dx 
 \\ \notag &
=\int_{\Sr \ts \S^{2} }\g _{s}(x,y)\cdot \g _{s}(x,z)\rho ^{3}(x,y,z)dxdydz + 
\int_{\Sr \ts \SS } |\g _{s}(x,y)|^{2} \rho ^{2}(x,y)dxdy 
.\end{align*}
Hence \eqref{:51a} implies $ \ggs \in L^2_{\mathrm{loc}}(\muone ) $. 
\bbbbb 

Let $ \rrNs (x,\mathsf{y})$ be as in \eqref{:44z}. 
Let $ \muN $ and $ \muNx $ be as in \sref{s:4} 
with $ n $-point correlation functions $ \rN $ and $ \rN _{x} $, respectively. 
\begin{lem} \label{l:52}
\eqref{:44e} with $ \hat{p} = 2 $ follows from the following: 
\begin{align}\label{:52a}&
\2  
\left| \mathrm{E}^{\muNx }
[\rrNs (x,\mathsf{y})] - 
\mathrm{E}^{\muN } [\rrNs (x,\mathsf{y})]\right|
= 0
,\\
\label{:52b}&
\2  
\left|\mathrm{Var}^{\muNx } [\rrNs (x,\mathsf{y})] -
\mathrm{Var}^{\muN } [\rrNs (x,\mathsf{y})] \right|= 0
,\\\label{:52c}&
\2  
\left| \mathrm{E}^{\muN } [\rrNs (x,\mathsf{y})] - \www  
\right| = 0
,\\
\label{:52d}&
\2  \mathrm{Var}^{\muN } [\rrNs (x,\mathsf{y})]= 0.
\end{align}
\end{lem}
\aaaaa 
By \eqref{:26e} we have 
\begin{align}\label{:52e}
 \int_{\Sr \ts \SSS } & |\rrNs (x,\mathsf{y})|^{2} d\muNone  = 
\int_{\Sq }\mathrm{E}^{\muNx }[|\rrNs (x,\mathsf{y})|^{2}]\, \rNone (x)\, dx 
\\ \notag &= \int_{\Sq } \{
| \mathrm{E}^{\muNx } [\rrNs (x,\mathsf{y})] |^{2} + 
\mathrm{Var}^{\muNx } [\rrNs (x,\mathsf{y})]\}  \, \rNone  (x)\, dx 
.\end{align}
By \eqref{:41b} we see that $\{\rNone \}_{N}$ is uniformly bounded. 
So \eqref{:44e} with $ \hat{p} = 2 $ follows from \eqref{:52a}--\eqref{:52e}. 
\bbbbb 

We give a sufficient condition of \eqref{:52a}--\eqref{:52d} in terms of correlation functions. 
\begin{lem} \label{l:53}
We set  $ \SS _{s\infty}^{x} = \{ y \in \SS ;s \le |x-y| < \infty \}$. Then 
\eqref{:52a}--\eqref{:52d} follow from \eqref{:53a}--\eqref{:53d}.  
\begin{align}\label{:53a}&
\2  |\int_{\SS _{s\infty}^{x}}\gN (x,y)
\{ \rNone _{x} (y)-\rNone (y)\}dy | = 0
,\\%
\label{:53b}&
\2  
|\int_{\SS _{s\infty}^{x}}|\gN (x,y)|^{2} \{ \rNone _{x} (y)-\rNone (y)\} dy
\notag \\ & \quad \quad
-\int_{(\SS _{s\infty}^{x})^{2} }
\gN (x,y)\cdot \gN (x,z) \{ \rNtwo _{x} (y,z)-\rNtwo (y,z)\}dydz |= 0 
,\\
\label{:53c}&
\2  
|\int_{\SS _{s\infty}^{x}}\{\vN (x,y) +\gN (x,y)\rNone (y)\}dy - \ww (x)| = 0
,\\\label{:53d}& 
\2  
|\int_{\SS _{s\infty}^{x}}|\gN (x,y)|^{2} \rNone (y)dy 
\\ \notag &\quad \quad \quad \quad \quad \quad \quad \quad \quad \quad \quad 
- \int_{(\SS _{s\infty}^{x})^{2}} 
\gN (x,y)\cdot \gN (x,z)\rNtwo (y,z) dydz |= 0 
.\end{align}
\end{lem}
\aaaaa 
This lemma is clear from the standard calculation of correlation functions 
combined with \eqref{:44z}. 
\bbbbb

\section{Log derivative of the Ginibre random point field. } \label{s:6}
In this section we calculate the log derivative 
$ \dgin $ of the Ginibre random point field $ \mug $. 
Let $ \mugN $ be the probability measure on $ \SSS $ whose $ n $-point correlation function $ \rgN $ is given by 
\begin{align}\label{:60a}&
\rgN (\mathbf{x}_n) = \det [\kgN (x_i,x_j)]_{1\le i,j \le n }
.\end{align}
Here $ \mathbf{x}_n=(x_1,\ldots,x_n)$ and $ \kgN $ is the kernel defined by 
\begin{align}\label{:60b}&
\kgN (x,y)= \frac{1}{\pi } e ^{-(|x|^{2} + |y|^{2})/2} \{ 
\sum_{n=0}^{N-1} \frac{(x\bar{y})^n}{n!}\} 
.\end{align}
We easily see that 
\begin{align}\label{:60c}&
|\kgN (x,y)| \le 
\frac{1}{\pi } e^{-||x|-|y||/2}
.\end{align}
By \eqref{:60a} and \eqref{:60b} the 1-point correlation function $ \rg ^{N,1} $ is given by 
\begin{align}\label{:60d}&
\rg ^{N,1} (x) = \frac{1}{\pi } e ^{-|x|^{2}}
\{ \sum_{k=0}^{N-1} \frac{|x|^{2k}}{k!} \} 
.\end{align}
Moreover, it holds that 
$ \rg ^{N,n}= 0 $ if $ n \ge N+1 $, and that for $ 2 \le n \le N $ 
\begin{align}\label{:60e}&
\rg ^{N,n}(\mathbf{x}_n) = \frac{1}{\pi ^n}\{\prod_{k=0}^{n-1} 
\frac{1}{k!}\}\, e ^{ -\sum_{k=1}^{n}|x_k|^{2} }\prod _{ i < j }^{n}|x_i-x_j|^{2} 
.\end{align}
Note that $ \mugN (\{\mathsf{s}(\SS )=N \})= 1 $. 
So by \eqref{:60e} with $ n=N$ we deduce that 
\begin{align}\label{:60f}&
\dginN (x,\mathsf{y})= -2x + \sum_{i=1}^{N-1}\frac{2(x-y_i)}{|x-y_i|^{2}} 
\quad \quad \Y 
.\end{align}
Let $ \mugone $ be the measure defined by \eqref{:26e} for $ \mug $. 
\begin{thm} \label{l:61} 
The log derivative $\dgin \in \Llocp (\mugone )$ exists for any $ 1 \le p < 2 $ and is given by 
\begin{align}\label{:61a}&
\dgin (x,\mathsf{y})= \limi{r} \sum_{\xyi < r }\frac{2(x-y_i)}{|x-y_i|^{2}}
\quad \quad (\mathsf{y}=\sum _i\delta_{y_i})
.\end{align}
The convergence of the series in the right-hand side takes place 
in $ \Llocp (\mugone )$. 
\end{thm}

To prove \tref{l:61} we use \tref{l:44}. 
So we check all the conditions in \tref{l:44}. 
For this purpose we first prepare several lemmas. 
\begin{lem} \label{l:62}
\eqref{:41a} and \eqref{:41b} hold. 
\eqref{:42c}, \eqref{:44a}, and \eqref{:44b} hold by taking 
\begin{align*}&
\uN (x)= \uu (x) = -2x , \quad  \vN (x,y) =\vvv (x,y) = 0 , \quad  \www = 2x 
,\\ &
\gN (x,y)= \g (x,y)= 2(x-y)/ |x-y|^{2} , \quad \hat{p} = 2 
.\end{align*}
\end{lem}
\aaaaa 
\eqref{:41a} follows immediately from \eqref{:21a}, \eqref{:21b}, \eqref{:60a}, and \eqref{:60b}. 
Let $ v_{i} $ be the norm of the $ i $th row vector of the matrix 
$ [\kgN (x_i,x_j)]_{1\le i,j \le n }$. Then by \eqref{:60c} we deduce that 
$ v_{i}\le \sqrt{n}/\pi $. So we deduce from \eqref{:60a} that 
$ \rgN (\mathbf{x}_n)\le \prod_{i=1}^{n} v_{i} \le (\sqrt{n}/\pi )^n $, 
which implies \eqref{:41b}.  
\eqref{:42c}, \eqref{:44a}, and \eqref{:44b} are trivial. 
\bbbbb 

By \lref{l:62} it only remains to prove \eqref{:44e} with $ \hat{p} = 2 $ 
for \tref{l:61}. 
By the argument in \sref{s:5} we see that \eqref{:44e} follows from 
\eqref{:52a}--\eqref{:52d}, which we will check below. 

It is known that the Palm measure conditioned at $ x $ of 
determinantal random point fields with kernel $ K $ 
is again a determinantal random point field with kernel 
$ K_x (y,z)= K(y,z) - \{K(y,x)K(x,z)/K(x,x)\}$ 
(see \cite[Theorem 1.7]{shirai-t}). 
Applying this to $ \mugNone $ we deduce that the kernel $ \kgNx $ 
of the Palm measure $ \mugNx $ is then given by 
\begin{align}\label{:63a}&
\kgNx (y,z)= \kgN (y,z)- \frac{\kgN (y,x)\kgN (x,z)}{\kgN (x,x)}
.\end{align}
Let $\Ct{;52e}=({1}/{\pi })
\sup _{x\in\Sq } e^{5|x|^{2}} $. 
Then by \eqref{:60b}, \eqref{:60c} and \eqref{:63a} 
we deduce that 
\begin{align}\label{:63b}
&|\kgNx (y,z) -\kgN (y,z)| = 
|\frac{\kgN (y,x)\kgN (x,z)}{\kgN (x,x)}|
\le \cref{;52e}e^{-(|y|^{2} + |z|^{2})/8}
.\end{align}
\begin{lem} \label{l:63}
\eqref{:52a} and \eqref{:52b} hold. 
\end{lem}
\aaaaa 
Since $ \mugN $ and $ \mugNx $ are determinantal random point fields with kernels $ \kgN $ and $ \kgNx $ respectively, 
their 1-point correlation functions $ \rgNone $ and $ \rgNonex $ are given by 
$  \rgNone (y) = \kgN (y,y)$ and $ \rgNonex (y) = \kgNx (y,y) $. Moreover, 
\begin{align*}&
\rgNtwo (y,z) = \kgN (y,y)\kgN (z,z)- \kgN (y,z)\kgN (z,y),\\  &
\rgNtwox (y,z) = \kgNx (y,y)\kgNx (z,z)- \kgNx (y,z)\kgNx (z,y)
.\end{align*}
Hence \eqref{:53a} and \eqref{:53b} follow from \eqref{:63b} 
and $ \gN (x,y)= {2(x-y)}/{|x-y|^{2}}$, which implies 
\eqref{:52a} and \eqref{:52b}. 
\bbbbb 

\begin{lem} \label{l:64}
Set $ \SS _{s\infty}=\{ s\le |y| < \infty \} $. Suppose $ \q < s $. Then 
\begin{align}\label{:64a}&
\int_{\SS _{s\infty}} \frac{2(x-y)}{|x-y|^{2}}
\rg ^{N,1}(y)dy =0 \quad \text{ for } x\in \Sq 
.\end{align}
\end{lem}
\aaaaa 
We regard $ x,y\in\R ^{2} $ as $ x,y \in \mathbb{C}$ 
and $ \bar{\cdot}$ denotes the complex conjugate. Then $ {(x-y)}/{|x-y|^{2}}={1}/({\bar{x}-\bar{y}})$. 
Recall that $ \rg ^{N,1}(y)=\rg ^{N,1}(|y|)$. 
Then 
\begin{align*}
\int_{\SS _{s\infty}}\frac{x-y}{|x-y|^{2}}
\rg ^{N,1}(y)dy &= 
\int_{\SS _{s\infty}}\frac{1}{\bar{x}-\bar{y}}
\rg ^{N,1}(|y|)dy 
\\ \notag &
= - \int_{\SS _{s\infty}}
\sum_{m=0}^{\infty} 
\bar{x}^{m}(\frac{1}{\bar{y}})^{m+1} \rg ^{N,1}(|y|)dy 
\quad \text{ by }\frac{|\bar{x}|}{|\bar{y}|}< 
\frac{\q }{s} <1 
\\ \notag &
= - \sum_{m=0}^{\infty} \bar{x}^{m}\int_{\SS _{s\infty}}
\frac{y^{m+1}}{|y|^{2(m+1)}} \rg ^{N,1}(|y|)dy
= 0
,\end{align*}
which implies \eqref{:64a}. Here we used $ 1/\bar{y}= y/|y|^{2} $ 
and \eqref{:60d} for the last line. 
\bbbbb 

Let $\SS _{s\infty}^{x}= \{ s \le |x-y| < \infty \}$. 
Note that $ \SS _{s\infty}=\SS _{s\infty}^{0}$. We set 
\begin{align}\label{:65a}&
T_s^x = \SS _{s\infty} \backslash \SS _{s\infty}^{x}, \quad 
U_s^x =\SS _{s\infty}^{x}\backslash \SS _{s\infty} 
.\end{align}
\begin{lem} \label{l:66}
\eqref{:52c} holds with $ \www = 2x $. 
\end{lem}
\aaaaa   
By \eqref{:44z}, $\vN (x,y) = 0 $, and $\gN (x,y)= 2(x-y)/ |x-y|^{2}$, we have 
\begin{align}\label{:66z}
\mathrm{E}^{\muN } [\rrNs (x,\mathsf{y})] = &
\int_{\SS _{s\infty}^{x}} \frac{2(x-y)}{|x-y|^{2}} \rg ^{N,1}(y)dy 
\\\notag 
=&-\int_{T_s^x}\frac{2(x-y)}{|x-y|^{2}}
\rg ^{N,1}(y)dy + 
\int_{U_s^x}\frac{2(x-y)}{|x-y|^{2}}
\rg ^{N,1}(y)dy 
\\ \notag {\to }&  
-\int_{T_s^x}\frac{2(x-y)}{|x-y|^{2}}\frac{1}{\pi }dy + 
\int_{U_s^x}\frac{2(x-y)}{|x-y|^{2}}\frac{1}{\pi }dy 
.\end{align}
uniformly in $ x\in \Sq $ as $ N\to \infty $. 
We used here \eqref{:64a} and \eqref{:65a} for the second line, and 
 \eqref{:60d} for the third one. 
By a direct calculation we obtain 
\begin{align}\label{:65b}&
\limi{s}\sup_{x\in\Sq }
|-\int_{T_s^x}\frac{2(x-y)}{|x-y|^{2}}\frac{1}{\pi }dy + 
\int_{U_s^x}\frac{2(x-y)}{|x-y|^{2}}\frac{1}{\pi }dy -2x | = 0 
.\end{align}
Combining \eqref{:66z} and \eqref{:65b} we obtain \eqref{:52c}. 
\bbbbb

\begin{lem} \label{l:69}   Let 
$\mathsf{h}_{rs}$ be the function on $ \SS \ts \SSS $ defined by 
\begin{align}\label{:67a}&
\mathsf{h}_{rs}(x,\mathsf{y})=\sum_{r\le \xyi < s } 
\frac{2(x-y_i)}{|x-y_i|^{2}}\,  \lceil |x-y_i| \rceil 
,\end{align}
where $ \lceil\cdot \rceil $ is the minimal integer 
greater than or equal to $\cdot $ and $\mathsf{y}=\sum_{i}\delta _{y_i}$. 
Then 
\begin{align}\label{:68a}&
\supN \sup_{x\in\Sq }\mathrm{Var}^{\mugN }
[\mathsf{h}_{rs}(x,\mathsf{y})] 
= O(s) \quad \text{ for each }r > 0 
.\end{align}
Here $ f(s)=O(s)$ means 
$ \limsup _{s\to\infty}|f(s)|/s < \infty $. 
\end{lem}
\aaaaa 
Let $ \SS _{rs} = \SS _{s}\backslash\SS _{r}$. Let $ h_{rs}(z)=1_{\SS _{rs}} (z)
{2z\lceil |z| \rceil }/{|z|^{2}}$. Then 
$ \mathsf{h}_{rs} (x,\mathsf{y})=\sum_i h_{rs}(x-y_i)$ by \eqref{:67a}. 
By a standard calculation of determinantal random point fields, we deduce that 
\begin{align}\label{:69b} 
\mathrm{Var}^{\mugN }[\mathsf{h}_{rs}(x,\mathsf{y})] & = 
-\int_{\SS _{rs} ^{2}}
(h_{rs}(x-y),h_{rs}(x-z))_{\R^{2}} |\kgN (y,z)|^{2} dydz 
\\ \notag & \quad \quad \quad \quad \quad \quad \quad \quad \quad + 
\int_{\SS _{rs} }|h_{rs}(x-y)|^{2} \kgN (y,y) dy 
.\end{align}

We set $\SS _{rs}^{x}= \{ r \le |x-y|< s \}$. 
By a direct calculation we have 
\begin{align}\label{:69c}&
1_{\SS _{rs}^{x}\cap\SS _{rs}}(y)| h_{rs} (x-y)-h_{rs}(-y)|  /2
\\ \notag =  &
1_{\SS _{rs}^{x}\cap\SS _{rs}}(y)\left| 
\frac{x-y}{|x-y|^{2}} \lceil |x-y| \rceil  + 
\frac{y}{|y|^{2}}\lceil |y| \rceil \right| 
\\ \notag = & 
1_{\SS _{rs}^{x}\cap\SS _{rs}} (y)\left| 
\{\frac{x-y}{|x-y|^{2} } + \frac{y}{|y|^{2}}\} \lceil |x-y| \rceil - 
\frac{y}{|y|^{2}} \{\lceil |x-y| \rceil - \lceil |y| \rceil \}\right|
\\ \notag \le & 
\cref{;54a} 1_{\SS _{rs}^{x}\cap\SS _{rs}} (y) /{|y|}
.\end{align}
Here $ \Ct{;54a}=\cref{;54a}(r)$ is the finite constant defined by 
\begin{align*}&
\cref{;54a}= 
\sup_{ x\in\Sq  \atop  y \in \SS _{r\infty }^{x}\cap\SS _{r\infty }}
 \left| 
\{\frac{x-y}{|x-y|^{2} } + \frac{y}{|y|^{2}} \}
\lceil |x-y| \rceil  - 
\frac{y}{|y|^{2}} \{\lceil |x-y| \rceil - \lceil |y| \rceil \}\right|
 |y|
.\end{align*}
Let $ \Ct{;69d}= \max\{2  \cref{;54a} + 2(r+1)/r\} $. 
Then by $  |h_{rs}(z)| \le 1_{\SS _{rs}} (z)\cdot 2(r+1)/r $ and 
\eqref{:69c}, 
we deduce that for all $ x \in \Sr $ and $ y \in \SS $
\begin{align}\label{:69e}&
| h_{rs} (x-y)-h_{rs}(-y)| \le \cref{;69d}
\{\frac{1_{\SS _{rs}^{x}\cap\SS _{rs}} (y)}{|y|} + 
1_{\SS _{rs}^{x}\backslash\SS _{rs}} (y) + 1_{\SS _{rs}\backslash \SS _{rs}^{x}} (y) \}
.\end{align}

By \eqref{:60c}, \eqref{:69b}, and \eqref{:69e} we easily deduce that 
\begin{align}\label{:69a}&
\sup _{N\in\N}\sup_{x\in\Sq }| 
\mathrm{Var}^{\mugN }[\mathsf{h}_{rs}(x,\mathsf{y})] - 
\mathrm{Var}^{\mugN }[\mathsf{h}_{rs}(0,\mathsf{y})] |= O(s)
.\end{align}
By applying \cite[Lemma 9.2]{o.rm} to $\mathsf{h}_{rs}(0,\mathsf{y}) $ we have 
\begin{align}\label{:69d}&
\supN \mathrm{Var}^{\mugN } [\mathsf{h}_{rs}(0,\mathsf{y})] = O(s)
.\end{align}
Hence \eqref{:68a} follows immediately from \eqref{:69a} and \eqref{:69d}. 
\bbbbb

Let $ \ggrs (x,\mathsf{y})= \sum_{r\le \xyi < s }{2(x-y_i)}/{|x-y_i|^{2}}
$. 
Then we easily deduce that 
\begin{align}\label{:67c}&
\ggrs  = \frac{\mathsf{h}_{rs}}{s} + 
\sum_{t=r+1}^{s-1} \frac{\mathsf{h}_{rt}}{t(t+1)}
.\end{align}
\begin{lem} \label{l:67}
$\mathsf{g}_{r\infty}(x,\cdot )= \limi{s}\ggrs (x,\cdot )$ exists in 
$ L^{2}(\mugN )$ for all $ x $ and 
\begin{align}\label{:67d}&
\mathsf{g}_{r\infty }(x,\cdot )
 = \sum_{t=r+1}^{\infty } \frac {\mathsf{h}_{rt}(x,\cdot )}{t(t+1)} 
\quad \text{ in } L^{2}(\mugN ) \quad \text{ for all }x
.\end{align}
\end{lem}
\aaaaa  
By $ \mugN (\mathsf{s}(\SS )= N) = 1 $ we see that 
$ \limi{s}{\mathsf{g}_{rs}}(x,\cdot )$ and 
$ \limi{s}{\mathsf{h}_{rs}}(x,\cdot )$ exist 
in $ L^{2}(\mugN )$ for all $ x $. 
Hence $ \limi{s}{\mathsf{h}_{rs}}/{s}= 0 $ in $ L^{2}(\mugN )$ for all $ x $. 
This together with \eqref{:67c} implies \eqref{:67d}. 
\bbbbb

\begin{lem} \label{l:68} 
\eqref{:52d} holds. 
\end{lem}
\aaaaa  
Let 
$ \hat{g}_{s\infty}= 
\3 \mathrm{Var}^{\mugN }
[\mathsf{g}_{s\infty }(x,\mathsf{y})]^{1/2} $. 
Let $\hat{h}_{rs}$ be defined similarly to 
$\hat{g}_{s\infty}$ by replacing 
$\mathsf{g}_{s\infty }$ by $ \mathsf{h}_{rs}$. 
By \eqref{:67a} we see that 
$ \mathsf{h}_{st}=-\mathsf{h}_{rs}+\mathsf{h}_{rt}$. 
So $ \hat{h}_{st}\le \hat{h}_{rs} + \hat{h}_{rt}$. 
Hence, by \eqref{:67d} we deduce that 
\begin{align}\label{68b}
\hat{g}_{s\infty} &\le 
\sum_{t=s+1}^{\infty } \frac{\hat{h}_{st}}{t(t+1)}
\le 
\sum_{t=s+1}^{\infty } \frac{\hat{h}_{rs}}{t(t+1)}+
\sum_{t=s+1}^{\infty } \frac{\hat{h}_{rt}}{t(t+1)}
.\end{align}
By \eqref{:68a} we deduce that $ \hat{h}_{rs}= O (\sqrt s)$. 
Combining this with \eqref{68b} we deduce that $ \limi{s} \hat{g}_{s\infty}=0 $, 
which yields \eqref{:52d}. 
\bbbbb 

\noindent {\em Proof of \tref{l:61}.} 
We use \tref{l:44} to prove \tref{l:61}. 
So we check that $ \{ \mugN \}$ satisfies the conditions in \tref{l:44}. 
By \lref{l:62} it only remains to prove \eqref{:44e}. 
Recall that \eqref{:44e} follows from \eqref{:52a}--\eqref{:52d}. 
We obtain \eqref{:52a} and \eqref{:52b} by \lref{l:63}. 
\eqref{:52c} follows from \lref{l:66}. 
\eqref{:52d} follows from \lref{l:68}. 
\qed

\section{Proof of Theorems \ref{l:21}--\ref{l:23}. }\label{s:66}
In this section we prove Theorems \ref{l:21}--\ref{l:23}. 
We recall that we took $ \vvv (x,y)=0$ in \lref{l:62}. So we set 
\begin{align} 
\label{:71z}&
\ggrs (x,\mathsf{y}) = \sum_{r \le |x-y_{i}|< s }{2(x-y_i)}/{|x-y_i|^{2}}
, \quad  \ggs =\mathsf{g}_{0s}
,\\\label{:71a}&
\ggtilders (x,\mathsf{y}) = \sum_{r \le |y_{i}|< s }{2(x-y_i)}/{|x-y_i|^{2}}
, \quad  \ggtildes =\ggtilde _{0s}
,\end{align}
where $ x \in \SS $ and $ \mathsf{y}=\sum_{i}\delta_{y_i}$. 

\begin{lem} \label{l:71}
\begin{align}\label{:71b}&
\limi{s}\{\ggs (x,\cdot ) - \ggtildes (x,\cdot ) \} = -2x 
\text{ in }L^{2}(\mug ) 
\text{ compact uniformly in $ x $.}
\end{align}
\end{lem}
\aaaaa  
Let $ T_s^x $ and $ U_s^x $ be as in \eqref{:65a}. 
By $ \rg ^{1}(x)= 1/\pi $ and \eqref{:65b} we deduce that 
\begin{align}\label{:71c}&
\limi{s} \mathrm{E}^{\mug }[
\ggs (x,\mathsf{y}) - \ggtildes (x,\mathsf{y})]
\\ \notag = & \limi{s} \{
\int_{T_s^x}\frac{2(x-y)}{|x-y|^{2}}\rg ^{1}(y)dy - 
\int_{U_s^x}\frac{2(x-y)}{|x-y|^{2}}\rg ^{1}(y)dy \} 
\\ \notag = & -2x \quad \text{ compact uniformly in }x 
.\end{align}
By a similar equality to \eqref{:69b} with $ \kg (y,y)= 1/\pi $ we deduce that 
\begin{align}\label{:71d}&
\limi{s}\mathrm{Var}^{\mug }
[\ggs (x,\mathsf{y}) - \ggtildes (x,\mathsf{y})]
\\ \notag \le & 
\limi{s} \{ \int_{T_s^x}\frac{4}{|x-y|^{2}}\frac{1}{\pi }dy + 
\int_{U_s^x}\frac{4}{|x-y|^{2}}\frac{1}{\pi }dy \} 
\\ \notag 
= & \ 0 \quad \text{ compact uniformly in }x 
.\end{align}
By \eqref{:71c} and \eqref{:71d} we obtain \eqref{:71b}.   
\bbbbb

We next prove the identity  of the form  
\begin{align}\label{:72a}&
\limi{s} \ggs (x,\mathsf{y}) = -2x + \limi{s}\ggtildes (x,\mathsf{y})
.\end{align}
\begin{lem}\label{l:72}
\thetag{1} For all $ x \in \Sr $, $ \ggrs (x,\mathsf{y})$ and 
$ \ggtilders (x,\mathsf{y})$ converge in $ L^{2}(\mug )$ as $ s\to \infty $ compact uniformly in $ x \in \Sr $. 
Moreover, \eqref{:72a} holds in the sense that 
\begin{align}\label{:72b}&
\ggr (x,\mathsf{y}) + \limi{s}\ggrs (x,\mathsf{y}) = 
-2x + \ggtilder (x,\mathsf{y})+ \limi{s}\ggtilders (x,\mathsf{y})
.\end{align}
\thetag{2} For all $ x $, $ \ggs (x,\mathsf{y})$ and 
$ \ggtildes (x,\mathsf{y})$ converge in $ L^{2}(\mugx )$ as $ s\to \infty $ compact uniformly in $ x $. 
Moreover, \eqref{:72a} holds. 
\\
\thetag{3} 
$ \ggs $ and $ \ggtildes $ converge 
in $ \Lloctwo (\mugone )$  as $ s\to \infty $ and \eqref{:72a} holds. 
\end{lem}
\begin{rem}\label{r:61}
Note that 
$ \ggr (x,\cdot ), \ggtilder (x,\cdot ) \not\in L^{2}(\mug )$ 
because of the singularity at $ x $. 
So the statement of \thetag{1} is weaker than the others. 
\end{rem}
\aaaaa 
By \cite[Theorem 1.3]{o-s} we see that 
$ \mathrm{Var}^{\mug } [\mathsf{h}_{rs} (x,\mathsf{y})] = O ( s ) $ 
compact uniformly in $ x \in \Sr $. 
Since $\mathrm{E}^{\mug }[\mathsf{h}_{rs} (x,\mathsf{y})] = 0 $ 
for all $ x \in \Sr $, we deduce that 
$\mathrm{E}^{\mug }[|\mathsf{h}_{rs} (x,\mathsf{y})|^{2}]= 
 \mathrm{Var}^{\mug }[\mathsf{h}_{rs} (x,\mathsf{y})] = O ( s ) $ 
compact uniformly in $ x \in \Sr $. 
Hence $ \limi{s} \mathsf{h}_{rs} /s = 0 $ in $ L^{2}(\mug )$ 
compact uniformly in $ x \in \Sr $. From this, combined with \eqref{:67c}, we deduce that $ \mathsf{g}_{r\infty}:=\limi{s}\ggrs $ 
converges in $ L^{2}(\mug )$ 
compact uniformly in $ x \in \Sr $. 
So by \eqref{:71b} we obtain \eqref{:72b}. We have thus proved \thetag{1}. 

By \eqref{:63b} and \eqref{:60a} and a similar representation of correlation functions of $ \mugNx $ we deduce that the first statement of \thetag{2} 
follows from that of \thetag{1}. Since $ \ggr , \ggtilder \in L^{2}(\mugx )$, 
the second follows from \eqref{:72b}. So we obtain \thetag{2}. 

\thetag{3} follows from \thetag{2} and the relation 
$ \mugone (A\ts B)= \int_{A}\mugx (B)\rg ^{1}(x)dx $ 
with $ \rg ^{1}(x)=1/\pi $. 
\bbbbb

\noindent 
{\em Proof of Theorems \ref{l:21} and \ref{l:23}.} 
We use Theorems \ref{l:26} and \ref{l:27} to prove \tref{l:21}. 
We take $ \mu = \mug $ and 
$ \mathsf{b}(x,\mathsf{y})= (1/2)\lim_s \ggs (x,\mathsf{y})$, 
where $ \ggs $ is same as \eqref{:72a}. Moreover, $\sigma (x,\mathsf{y})$ is the unit matrix for all $(x,\mathsf{y})$. Hence $ \aaa  = \sigma ^{2} $ is also the unit matrix. 
We check that $ \mug $ satisfies \thetag{A.1}--\thetag{A.5} for these 
$ \sigma $ and $\mathsf{b}$. 

\thetag{A.1} and \thetag{A.5} are clear from \eqref{:21a} and \eqref{:21b}. 
\thetag{A.2} follows from \tref{l:61} and \lref{l:72} \thetag{3}. 
In \cite[Theorem 2.6]{o.rm} we proved that the closability in \thetag{A.3} holds for $ k=0 $. Indeed, we proved that $\mug $ is a quasi-Gibbs measure in the sense of \cite[Definition 2.1]{o.rm} and deduced the closability for $ k=0 $ from this. The closability for general $ k \in \N $ also follows in a similar fashion from the quasi-Gibbs property of $ \mug $. Since the kernel $ \kg $ is locally Lipschitz continuous, \thetag{A.4} immediately follows from \cite[Theorem 2.1]{o.col}. 

We thus see that $ \mug $ satisfies \thetag{A.1}--\thetag{A.5}. 
Hence Theorems \ref{l:21} and  \ref{l:23} follow from 
Theorems \ref{l:26} and \ref{l:27}, respectively. 
\qed 

\noindent 
{\em Proof of \tref{l:22}.} 
By \lref{l:72} \thetag{3} we see that \eqref{:72a} holds 
in $ \Lloctwo (\mugone )$. Hence we deduce that 
$ \bbb (x,\mathsf{y}) = -x + \tilde{\bbb }(x,\mathsf{y})$ 
in $ \Lloctwo (\mugone )$. This combined with \tref{l:21} implies 
\tref{l:22}. 
\qed

\section{Proof of Theorems \ref{l:24} and \ref{l:25}.}\label{s:8} 
In this section we prove Theorems \ref{l:24} and \ref{l:25} by using 
Theorems \ref{l:26} and \ref{l:27}. So we take $ \mu = \mub $ 
and prove that $ \mub $ satisfies \thetag{A.1}--\thetag{A.5}. 

\begin{lem} \label{l:81}
$ \mub $ $( \beta = 1,2,4 )$ satisfy 
\thetag{A.1}, \thetag{A.3}, \thetag{A.4}, and \thetag{A.5}.
\end{lem}
\aaaaa 
Since the correlation functions $\{\rb \}$ of $\mub $ have the expression \eqref{:24} and the kernels $ \Ksinb $ are bounded, \thetag{A.1} and \thetag{A.5} are clear. 

In \cite[Theorem 2.5]{o.rm} we proved that the closability in \thetag{A.3} holds for $ k=0 $. Indeed, we proved that $ \mub $ is a quasi-Gibbs measure and deduced the closability for $ k=0 $ from this. The closability for general $ k \ge 1 $ also follows from the quasi-Gibbs property of $ \mub $ in a similar fashion. Since the kernel $ \Ksinb $ is locally Lipschitz continuous, 
\thetag{A.4} follows from \cite[Theorem 2.1]{o.col}. 
\bbbbb 

By \lref{l:81} it only remains to prove \thetag{A.2}. 
Define $ \KsinNb (x) $ by \eqref{:91s}--\eqref{:91t} 
with the replacement of $ S(x)$ by 
$\SN (x)={\sin (\pi x)}/\{\nN {\sin (\pi x/\nN )}\}$. 
We set $ \RN = (-\nN /2,\nN /2] $ and 
\begin{align*}&
\KsinNb (x,y) = 1_{\RN }(x) \KsinNb (x-y)1_{\RN }(y) 
.\end{align*}
We take $ \muN $ in \thetag{A.2} to be the probability measure $\mubN $
on $ \SSS $ whose $ n $-point correlation function $ \rbN $ is given by 
\begin{align}\label{:82a} 
\rbN (\mathbf{x}) = \det [\KsinNb (x_i,x_j)]_{1 \le i,j \le n}
,\end{align}
where $ \mathbf{x}=(x_i)$. It is well known \cite{mehta} that 
$ \mubN (\mathsf{s}(\R )=N)=1 $ and that 
\begin{align}\label{:82b}&
\rbNN (\mathbf{x}) = 
\mathrm{const.} 
\prod _{i, j = 1, i < j }^{\nN } 1_{\RN }(x_i) 
|e^{2\pi \iii  x_i/\nN }- 
 e^{2\pi \iii  x_j/\nN }| ^{\beta } 1_{\RN }(x_j) 
.\end{align}
We can regard $ \RN $ as a torus and $ \mubN $ to be a translation invariant probability measure on the configuration space on the torus $ \RN $. 
The image measure of $ \mubN $ under the map 
$ \omega _{\nN }(x) = e^{2\pi \iii  x/\nN }$ gives the distributions of the eigenvalues of the random matrices called circular ensembles \cite{mehta}. 
We can rewrite \eqref{:82b} as 
\begin{align}\label{:82z}&
\rbNN (\mathbf{x}) = 
\mathrm{const.} 
\prod _{i, j = 1, i < j }^{\nN } 1_{\RN }(x_i) 
|\omega _{\nN }(x_i)- \omega _{\nN }(x_j)| ^{\beta } 1_{\RN }(x_j) 
.\end{align}

Taking \eqref{:82b} into consideration we set 
\begin{align}\label{:82d}
\gN (x,y)&= \PD{}{x} 
\log |e^{2\pi \iii  x/\nN }- e^{2\pi \iii  y/\nN }|^{\beta }
&&\text{ if } x,y\in (-\nN /2, \nN /2) 
\\ \notag & = 0 &&\text{ otherwise}
.\end{align}
Then we can easily check that 
\begin{align}\label{:82v}&
\limsupi{N}
\sup_{x\in\Sq } |\int_{\SS _{s\infty}^{x}}\gN (x,y) dy | = o(s) 
\quad (s\to\infty )
\end{align}
and that there exists a constant $ \Ct{;72g}$ such that 
\begin{align}\label{:82g}&
\sup_{\nN \ge 8r }\sup_{ x \in \Sq } 
|\gN (x,y) | \le \cref{;72g}\min\{1,1 /{|y|}\} 
\quad \text{ for all } |y| > 2r 
.\end{align}

\begin{thm} \label{l:82}
Suppose $ \beta = 1,2,4$. Then the log derivative $\dlog ^{\mub }$ exists 
in $\Llocp (\mub ^{1})$ for any $ 1< p < 2$. 
Moreover $\dlog ^{\mub }$ is given by 
\begin{align}\label{:82c}&
\dlog ^{\mub }(x,\mathsf{y})= \limi{r} \sum_{|x-y_i|< r }\frac{\beta }{x-y_i}
\quad \quad (\mathsf{y}=\sum_i\delta _{y_i})
.\end{align}
\end{thm}
\aaaaa  
We use \tref{l:44} to prove \tref{l:82}. 
So we check the conditions of \tref{l:44}. We take 
$ \uu (x)=\www = 0 $, $\uN (x)= \delta_{-N/2}(x)-\delta_{N/2}(x)$, 
where $ \delta_{\pm N/2}(x)$ are delta measures, and 
$\vN (x,y)= \vvv (x,y)=  0 $.  
We set $ \gN $ as \eqref{:82d} and $\g (x,y)= 2/(x-y)$.

The conditions \eqref{:41a} and \eqref{:41b} follow from \eqref{:82a} and 
the definition of $ \KsinNb $. \eqref{:42c} and \eqref{:44a} are clear. 
For $ \beta = 2 , 4$, the condition \eqref{:44b} with 
$ \hat{p} = 2 $ follows from \eqref{:24}, $\g (x,y)= 2/(x-y)$, 
\eqref{:82d}, and \lref{l:51}. For $ \beta = 1 $ one can check that 
\eqref{:44b} with $ 1<\hat{p}<2 $ holds by the H\"{o}lder inequality in addition to the above.

We next prove \eqref{:44e}. 
For this it is sufficient to check \eqref{:53a}--\eqref{:53d} by \lref{l:53}. 
Let $ \mubNx $ be the Palm measure of $ \mubN $ conditioned at $ x \in \RN $ 
and let $ \rbNxn $ be its $ n $-point correlation function. 
Then $ \mubNx $ has a  determinantal structure with kernel 
\begin{align}\label{:82k}&
\KsinNbx (y,z)= \KsinNb (y,z) - \KsinNb (y,x)\KsinNb (x,z)/ \KsinNb (x,x)
.\end{align}
When $ \beta = 2$, \eqref{:82k} follows from \cite[Theorem 1.7]{shirai-t}. 
When $ \beta = 1,4 $, one can also check \eqref{:82k}. 
By \eqref{:91q} and \eqref{:91s}--\eqref{:91t} 
we easily see that $\KsinNb (x,x)=1_{\RN }(x) $. 
Hence \eqref{:82k} implies that for $ x\in \RN $ and $ y,z \in \R $ 
\begin{align}\label{:82n}&
\KsinNbx (y,z)= \KsinNb (y,z) - \KsinNb (y,x)\KsinNb (x,z)
.\end{align}

By \eqref{:82a} and \eqref{:82n} 
we see that for $ x\in \RN $ and $ y,z \in \R $ 
\begin{align}\label{:82u}&
\rbNxone (y)- \rbNone (y) = 
- [\KsinNb (y,x)\KsinNb (x,y)]^{(0)}
.\end{align}
Here $ [\cdot ]^{(0)}$ means the scaler part of quaternions $ \cdot $ in the sense of the Appendix. When $\beta = 2 $, $ [\cdot]^{(0)}=\cdot $ because $ \cdot $ are complex numbers. 
By \eqref{:82u}, \eqref{:91x} and \eqref{:91y} 
there exists a constant $ \Ct{;72e}$ satisfying 
\begin{align}
\label{:82e}
\sup_{\nN \ge 8r}\sup_{ x \in \Sq }|\rbNxone (y)-&  \rbNone (y)| 
\le 
\cref{;72e}\min \{1,{1}/{|y|} \}
\quad \text{ for all } |y| >2r 
.\end{align}
By \eqref{:82g} and \eqref{:82e} we obtain \eqref{:53a}. 

By \eqref{:82a} and \eqref{:82n} 
we see that for $ x\in \RN $ and $ y,z \in \R $ 
\begin{align}\label{:82p}& 
\rbNxtwo (y,z)- \rbNtwo (y,z)
\\ \notag 
= &
-[\KsinNb (y,x)\KsinNb (x,y)]^{(0)}- [\KsinNb (z,x)\KsinNb (x,z)]^{(0)}
\\ \notag &
+[\KsinNb  (y,x)\KsinNb (x,z)\KsinNb (z,y) ]^{(0)}
+[\KsinNb (y,z)\KsinNb (z,x)\KsinNb (x,y)]^{(0)}
.\end{align}
Then by \eqref{:82v}, \eqref{:82g}, \eqref{:91x} and \eqref{:91y} we see that 
as $ s \to \infty $
\begin{align}\notag &
\limsupi{N}\sup_{x\in\Sq } \left|\int_{(\SS _{s\infty}^{x})^{2}}
 \gN (x,y)\gN (x,z) [\KsinNb (y,x)\KsinNb (x,y)]^{(0)}dydz \right| 
\\ \notag  \le &
\ o(s) \cdot \limsupi{N} 
\sup_{x\in\Sq } \int_{\SS _{s\infty}^{x}}\left|
 \gN (x,y)[\KsinNb (y,x)\KsinNb (x,y)]^{(0)} \right| dy  
 \quad \text{ by }\eqref{:82v} 
\\ \label{:82q} 
= & \ o(s) \quad \quad \quad \quad \quad \quad \quad 
 \quad \text{ by }\eqref{:82g},\ \eqref{:91x}, \text{ and } \eqref{:91y}   
.\end{align}
By \eqref{:82g}, \eqref{:91x} and \eqref{:91z} we deduce that 
\begin{align}\notag&
\limsupi{N}\sup_{x\in\Sq }
\int_{(\SS _{s\infty}^{x})^{2}}| \gN (x,y)\gN (x,z) 
[\KsinNb (y,x)\KsinNb (x,z)\KsinNb (z,y)]^{(0)}|dydz 
\\\label{:82r}&= 
 o(s) \quad ( s \to \infty ) 
.\end{align}
We therefore obtain \eqref{:53b} from 
\eqref{:82p}, \eqref{:82q} and \eqref{:82r}. 

By $ \rbNone  = 1_{\RN }$ and \eqref{:82v} we obtain  
\eqref{:53c} because $ \vN (x,y) = \www = 0 $.

By \eqref{:82a} and $ \KsinNb (y,y) = 1_{\RN }(y)$ we deduce that 
\begin{align}\label{:82j}&
\rbNtwo (y,z)= 1_{\RN } (y)1_{\RN } (z)
\{ 1  - [\KsinNb (y,z)\KsinNb (z,y)]^{(0)}\} 
.\end{align}
So we deduce from 
\eqref{:82v}, \eqref{:91x}, \eqref{:91y}, and \eqref{:82j} that 
\begin{align}\label{:82l}&
\4  
|\int_{(\SS _{s\infty}^{x})^{2} }\gN (x,y)\gN (x,z)\rbNtwo (y,z) dydz | 
\\ \notag 
= & \ o(s)+ \4  
|\int_{(\SS _{s\infty}^{x})^{2}}\gN (x,y)\gN (x,z)
[\KsinNb (y,z)\KsinNb (z,y)]^{(0)} dydz |
\\ \notag 
= &\ o(s) \quad ( s \to \infty )
.\end{align}
Since $ \rbNone (y)= 1_{\RN }(y)$, we deduce from \eqref{:82g} that 
\begin{align}\label{:82m}&
\4 
\int_{\SS _{s\infty}^{x}}|\gN (x,y)|^{2} \rbNone (y) dy = 
\ o(s) \quad ( s \to \infty )
. \end{align}
Hence by \eqref{:82l} and \eqref{:82m} we obtain \eqref{:53d}. 
\bbbbb 

\begin{lem} \label{l:83}
Suppose $ \beta = 2,4$. 
Then $ \dlog ^{\mub }\in \Lloctwo (\mub ^{1})$. 
\end{lem}
\aaaaa 
Let $ \mathsf{g}_{s}(x,\mathsf{y})= \sum_{|x-y_i|< s }{2}/{(x-y_i)}$. 
Then by \tref{l:82} it is sufficient for \lref{l:83} 
to prove $ \mathsf{g}_{s} $ 
converge in $ \Lloctwo (\mub ^{1}) $. 
Let $ \mu _{x}$ be the Palm measure of $ \mub $ conditioned at $ x $. 
Then since $ \mub $ are translation invariant, it is enough to show that 
$ \mathsf{g}_{s}(x,\mathsf{y}) $ converge in 
$ L^2(\mu _{x}) $ for each $ x $. 

Let $ \mathsf{h}_{s}(x,\mathsf{y})=\sum_{ \xyi < s } 
{2}\lceil |x-y_i| \rceil /(x-y_i) $. Then we see that 
\begin{align}\label{:83a}&
\ggs  = \frac{\mathsf{h}_{s}}{s} + 
\sum_{t=1}^{s-1} \frac{\mathsf{h}_{t}}{t(t+1)}
.\end{align}
By the calculation based on the 1 and 2-point correlation functions 
we can check $ \mathrm{E}^{\mu _{x}}[|\mathsf{h}_{s}|^2]\sim O(s)$. 
This combined with \eqref{:83a} completes the proof. 
\bbbbb

\noindent 
{\em Proof of Theorems \ref{l:24} and \ref{l:25}.} 
By \lref{l:81}, \tref{l:82}, and \lref{l:83} we see that 
$ \mub $ ($ \beta = 2,4$) satisfy \thetag{A.1}--\thetag{A.5}. 
Hence Theorems \ref{l:24} and \ref{l:25} follow from 
Theorems \ref{l:26} and \ref{l:27}, respectively. 
When $ \beta = 1 $, 
$ \dlog ^{\mub }\in \Llocp (\mub ^{1})$ for any $1< p < 2$ and 
$ \dlog ^{\mub }\not\in \Lloctwo (\mub ^{1})$. 
In this case we can justify \eqref{:26u} by using the localization, and 
we still have Theorems \ref{l:24} and \ref{l:25}. 
\qed

\section{Appendix. }\label{s:9} 
We begin by defining $ \Ksinb $ for $ \beta = 1,4 $. 
For this purpose, we recall the standard quaternion notation 
for $ 2\ts 2 $ matrices (see \cite[Ch.\ 2.4]{mehta}), 
\begin{align} \label{:91l}&
\mathbf{1} = 
\begin{bmatrix}
1& 0 \\ 0&1
\end{bmatrix},\quad 
\mathbf{e}_1 = 
\begin{bmatrix}
\iii & 0 \\ 0&-\iii 
\end{bmatrix},\quad 
\mathbf{e}_2 = 
\begin{bmatrix}
0&1\\-1&0
\end{bmatrix},\quad 
\mathbf{e}_3 = 
\begin{bmatrix}
0&\iii \\\iii &0
\end{bmatrix}
.\end{align}

A quaternion $ q $ is represented by 
$ q = q^{(0)}\mathbf{1} + q^{(1)}\mathbf{e}_1 + 
 q^{(2)}\mathbf{e}_2 + q^{(3)}\mathbf{e}_3 $, 
where $ q^{(i)} $ are complex numbers. There is a natural identification between the $ 2\ts 2 $ complex matrices and the quaternions given by 
\begin{align} \label{:91q}&
\begin{bmatrix}
a&b \\ c&d
\end{bmatrix} = 
\frac{1}{2}(a+d)\mathbf{1} -
\frac{\iii }{2}(a-d)\mathbf{e}_1 +
\frac{1}{2}(b-c)\mathbf{e}_2 -
\frac{\iii }{2}(b+c)\mathbf{e}_3 
.\end{align}
We denote by 
$ \Theta (\begin{bmatrix}a&b \\ c&d \end{bmatrix}) $ 
the quaternion defined by the right hand side of \eqref{:91q}. 

For a quaternion 
$ q = q^{(0)}\mathbf{1} +q^{(1)}\mathbf{e}_1 +
q^{(2)}\mathbf{e}_2 +q^{(3)}\mathbf{e}_3 $, 
we call $ q^{(0)} $ the scalar part of $ q $. 
A quaternion is called scalar if $ q^{(i)} = 0 $ for 
$ i = 1,2,3 $. We often identify a scalar quaternion 
$ q = q^{(0)}\mathbf{1}$ with the complex number $ q^{(0)} $. 

Let 
$ \bar{q} = q^{(0)}\mathbf{1}- 
\{ q^{(1)}\mathbf{e}_1 +q^{(2)}\mathbf{e}_2 +q^{(3)}\mathbf{e}_3 \}.
$ 
A quaternion matrix $ A = [a_{ij}]$ is called self-dual 
if $ a_{ij} = \bar{a}_{ji} $ for all $ i,j $. 
For a self-dual $ n\ts n $ quaternion matrix $ A = [a_{ij}]$ we set 
\begin{align} \label{:91p}&
\det A = \sum_{\sigma \in \mathfrak{S}_n }
\mathrm{sign} [\sigma ] 
\prod _{i = 1}^{L(\sigma )}
[a_{\sigma_i(1)\sigma_i(2)}a_{\sigma_i(2)\sigma_i(3)}
\cdots a_{\sigma_i(\ell -1)\sigma_i(\ell )}]^{(0)} 
.\end{align}
Here $ \sigma = \sigma_1\cdots\sigma_{L(\sigma )}$ 
is a decomposition of $ \sigma $ to products of 
the cyclic permutations $ \{\sigma_i\} $ 
with disjoint indices. We write 
$ \sigma_i = (\sigma_i(1),\sigma_i(2),\ldots,
\sigma_i(\ell )) $, where $ \ell $ 
is the length of the cyclic permutation $ \sigma_i $. 
The decomposition is unique up to the order of 
$ \{ \sigma _i \} $. 
As before $ [ \cdot ]^{(0)} $ means 
the scalar part of the quaternion $ \cdot $. 
It is known that the right hand side is well defined (see  
\cite[Section 5.1]{mehta}). 

We are now ready to introduce $ \Ksinb $. 
Let $ S(x) = \sin (\pi x)/ \pi x $ and define 
\begin{align} \label{:91s}&
\Ksinx (x) = \ \Theta (
\begin{bmatrix}
S(x)& \frac{d S }{d x}(x)
\\
\int_0^x S(y)dy -\frac{1}{2} \mathrm{sgn}(x)
& S(x)
\end{bmatrix})
,\\\label{:91u}&
\Ksiny (x) = S(x)
,\\ \label{:91t}& \Ksinz (x) = \Theta ( 
\begin{bmatrix}
S(2x)& \frac{d S }{d x}(2x) 
\\ 
\int_0^{2x} S(y)dy & S(2x)
\end{bmatrix})
.\end{align}
We thus clarify the meaning of \eqref{:24}. 

We set the kernels $ \KsinNb $ by \eqref{:91s}--\eqref{:91t} with the replacement of $ S(x) $ by 
$\SN (x)={\sin (\pi x)}/\{ \nN {\sin (\pi x/\nN )}\}$. 
Let $ \omega _{\nN }(x) =e^{2\pi \iii  x/\nN }$ as before, and set 
\begin{align}\label{:91x}&
\eta _{N}^{x,y} =1_{\RN }(x)1_{\RN }(y)
 \min \{ 1 , {1}/|\omega _{\nN }(x)-\omega _{\nN }(y)|\}
.\end{align}
Then by \eqref{:91q} and \eqref{:91s}--\eqref{:91t} 
there exist constants $\Ct{;728}$ and $ \Ct{;73e}$ such that 
\begin{align}\label{:91y}&
|[\KsinNb (x,y) \KsinNb (y,x)]^{(0)} |
 \le \cref{;728} 
\eta _{N}^{x,y}
,\\ \label{:91z}&
|[\KsinNb (x,y)\KsinNb (y,z)\KsinNb (z,x)]^{(0)}| \le 
\cref{;73e} \{ \eta _{N}^{x,y}\eta _{N}^{y,z} + 
\eta _{N}^{y,z}\eta _{N}^{z,x}  
+ \eta _{N}^{z,x}\eta _{N}^{x,y}  \} 
.\end{align}




%
%

\end{document}